\documentclass[3p]{elsarticle}

\usepackage{amssymb}
\usepackage{amsmath}
\usepackage{amsfonts}

\usepackage{booktabs}
\usepackage{multirow}
\usepackage{comment}

\usepackage{bbm}        % nicer symbols for sets
\usepackage[utf8x]{inputenc}
\usepackage{graphicx}

\newtheorem{remark}{Remark}
\usepackage[abs]{overpic}
\usepackage{todonotes}

\journal{Comp. Meth. Appl. Mech. Eng.}

\newcommand{\divergence}{\mbox{div}}

\newcommand{\normal}{n}

\newcommand{\bilinearform}[1]{\mathcal{#1}_h}
\newcommand{\operator}[1]{\mathbf{#1}}

\newcommand{\sumoverallelements}{\sum_{T \in \mathcal{T}_h}}

\newcommand{\jumpleft}{[\![}
\newcommand{\jumpright}{]\!]}

\newcommand{\stab}{\alpha}
\newcommand{\brokenHnormleft}{|\!|\!|}
\newcommand{\brokenHnormright}{|\!|\!|_{1,*}}
\newcommand{\brokenHnormIIleft}{|\!|\!|}
\newcommand{\brokenHnormIIright}{|\!|\!|_1}
\newcommand{\proj}{\Pi}

\newcommand{\rr}{\mathbb{R}}

\newcommand{\DGSymb}{V}
\newcommand{\HdivSymb}{W}
\newcommand{\FacetSymb}{F}
\newcommand{\PressureSymb}{Q}
\newcommand{\VelSymb}{U}

\newcommand{\DGVar}{w}
\newcommand{\HdivVar}{u_T}
\newcommand{\FacetVar}{u_F}
\newcommand{\FacetVarLO}{\bar{u}_F}
\newcommand{\FacetVarHO}{\lambda}
\newcommand{\FacetVarT}{\lambda_T}
\newcommand{\FacetVarTT}{\lambda_{T'}}
\newcommand{\PressureVar}{p}
\newcommand{\VelVar}{u}

\newcommand{\DGVec}{\mathbf{w}}

\newcommand{\PressureVec}{\mathbf{p}}
\newcommand{\VelVec}{\mathbf{u}}

\newcommand{\VelVecExpl}{\mathbf{v}}

\newcommand{\Test}{v}
\newcommand{\DGTest}{z}
\newcommand{\HdivTest}{v_T}
\newcommand{\FacetTest}{v_F}
\newcommand{\FacetTestLO}{\bar{v}_F}

\newcommand{\FacetTestT}{\mu_T}
\newcommand{\FacetTestTT}{\mu_{T'}}
\newcommand{\PressureTest}{q}
\newcommand{\VelTest}{v}

\newcommand{\DGSpace}{\DGSymb_h}
\newcommand{\HdivSpace}{\HdivSymb_h}
\newcommand{\FacetSpace}{\FacetSymb_h}
\newcommand{\FacetSpaceLO}{\overline{\FacetSymb}_h}
\newcommand{\PressureSpace}{\PressureSymb_h}
\newcommand{\VelSpace}{\VelSymb_h}

\newcommand{\transferop}{\operator{I}}
\newcommand{\transferopI}{\operator{I}^{T}}

\newcommand{\glowa}{n_0}
\newcommand{\glowb}{n_1}
\newcommand{\glowc}{n_2}
\newcommand{\glowd}{n_3}
\newcommand{\markdigits}[1]{#1}

\begin{document}

\begin{frontmatter}

% \runningheads{C.~Lehrenfeld and J.~Sch\"oberl}{Exactly divergence-free HDG methods for incompressible flows}

\title{High order exactly divergence-free Hybrid Discontinuous Galerkin Methods for unsteady incompressible flows}

\author[add1]{Christoph Lehrenfeld}
\ead{christoph.lehrenfeld@uni-muenster.de}
\author[add2]{Joachim Sch\"oberl}
\ead{joachim.schoeberl@tu-wien.ac.at}
\address[add1]{Institute for Computational and Applied Mathematics, WWU M\"unster, M\"unster, Germany}
\address[add2]{Institute for Analysis and Scientific Computing, TU Wien, Wien, Austria}
% \author{Christoph Lehrenfeld~and Joachim Sch\"oberl}

% \address{Institute for Computational and Applied Mathematics, WWU M\"unster, M\"unster, Germany}
% \corraddr{Christoph Lehrenfeld,
% Institute for Computational and Applied Mathematics at the University of M\"unster, Einsteinstrasse 62, D-48149 M\"unster, Germany.
%   E-mail: christoph.lehrenfeld@uni-muenster.de}

% \address{
%   \affilnum{1} Institute for Analysis and Scientific Computing, TU Wien,
%   Wiedner Hauptstrasse 8-10,
%   1040 Wien, Austria
% }

% \corraddr{Christoph Lehrenfeld,
%   Institute for Analysis and Scientific Computing, TU Wien,
%   Wiedner Hauptstrasse 8-10,
%   1040 Wien, Austria.
%   E-mail: christoph.lehrenfeld@gmail.com}

\begin{abstract}
 In this paper we present an efficient discretization method for the solution of the unsteady incompressible Navier-Stokes equations based on a high order (Hybrid) Discontinuous Galerkin formulation. The crucial component for the efficiency of the discretization method is the disctinction between stiff linear parts and less stiff non-linear parts with respect to their \emph{temporal and spatial} treatment.

  Exploiting the flexibility of \emph{operator-splitting time integration schemes} we combine two spatial discretizations which are tailored for two simpler sub-problems: a corresponding hyperbolic transport problem and an unsteady Stokes problem.

  For the hyperbolic transport problem a spatial discretization with an Upwind Discontinuous Galerkin method and an explicit treatment in the time integration scheme is rather natural and allows for an efficient implementation.
  The treatment of the Stokes part involves the solution of linear systems. In this case a discretization with Hybrid Discontinuous Galerkin methods is better suited. We consider such a discretization for the Stokes part with two important features: $H(\divergence)$-conforming finite elements to garantuee exactly divergence-free velocity solutions and a projection operator which reduces the number of globally coupled unknowns.
  We present the method, discuss implementational aspects and demonstrate the performance on two and three dimensional benchmark problems.
\end{abstract}

% \keywords{

% \maketitle

\end{frontmatter}

% \maketitle
% \tableofcontents

\section{Introduction}
\subsection{Problem statement}
We consider the numerical solution of the unsteady incompressible Navier Stokes equations in a velocity-pressure formulation:
\begin{equation} \label{eq:navierstokes}
\left\{
\begin{array}{r l @{\hspace*{0.1cm}} c l l}
 \frac{\partial}{\partial t} u +  \divergence ( - \nu \nabla u + u \! \otimes \! u  + p {I} ) 
 &=& f \! &\mbox{ in } \Omega \\
   \divergence \ \! u  &=& 0 &\mbox{ in } \Omega 
%   \\
%  				u &=&  u_D \! &\mbox{ on }
%\Gamma_D \\
%  				(\nu \nabla u  - p {I})
%\cdot \normal &=& 0 &\mbox{ on } \Gamma_{out}
\end{array} \right.
\end{equation}
% where $\Omega$ is a bounded domain in $\mathbb{R}^d$, $d=1,2,3$ with Lipschitz
% boundary $\partial \Omega = \Gamma_D \cup \Gamma_{out} $. The boundary
% parts $\Gamma_D$ and $\Gamma_{out}$ denote the part of the boundary
% where Dirichlet or natural outflow boundary conditions are prescribed
% respectively. $[t_0,T]\subset \mathbb{R}$ is the time interval , where $t_0$ is
% the start time where the initial condition $u_0$ is prescribed and $T$ is some
% later time. \\
with boundary conditions $u=u_D$ on $\Gamma_D \subset \partial \Omega$ and $(\nu \nabla u  - p {I}) \cdot \normal = 0$ on $\Gamma_{out} = \partial \Omega \setminus \Gamma_D$. Here, $\nu = const$ is the kinematic viscosity, ${u}$ the velocity, $p$ the pressure, and $f$ is an external body force. 
We consider a discretization to \eqref{eq:navierstokes} with high order (Hybrid) Discontinuous Galerkin (DG) finite element methods for complex geometries with underlying meshes consisting of possibly curved tetrahedra, hexahedra, prisms and pyramids. 

%
%% The incompressible Navier-Stokes equations contain multiple challenging aspects that have to be addressed when designing a discretization method.

%% \begin{itemize}
%% \item incompressibility constraint. (divergence-free ansatz functions typically not feasible). Saddle-point structure (LBB)
%% \item 
%% \end{itemize}
\subsection{Literature}
Since its introduction in the paper of Reed and Hill \cite{reed1973triangularmesh}, DG methods have been developed and used for hyperbolic problems \cite{lesaint1974finite,johnson1986analysis,bassi1997high,bassirebay97b,cockburn1998runge} and later extended to second-order elliptic (and parabolic) problems \cite{arnold2002unified,houston2002discontinuous,riviere2008discontinuous, dipietroern2012}.
DG methods are specifically popular for flow problems \cite{hesthaven2007nodal,karniadakis2013spectral}. For incompressible flows different finite element methods have been discussed in the literature 
\cite{girault2012finite,elman2014finite,donea2003finite}. Among these methods several DG formulations have been considered, e.g. \cite{toselli2002hp,schotzau2002mixed, girault2005discontinuous, cockburn2005locally, cockburn2007note}.

DG methods provide a flexibility which can be utilized for different purposes. In our case, the motivation to consider this type of discretizations is twofold. First, relevant flow problems that can be modeled with the incompressible Navier-Stokes equations are often convection dominated. Using the Upwind mechanism DG methods offer a natural way to devise stable discretizations of (dominating) convection.
Secondly, abandoning $H^1$-conforming finite element spaces, DG methods allow to consider (only) $H(\divergence)$-conforming finite elements (see \cite{brezzi2012mixed}) for (Navier-)Stokes problems. This is attractive as it facilitates the design of a discretization with exactly solenoidal solution which in turn implies energy-stability of the discretization for the Navier-Stokes problem, cf. \cite{cockburn2005locally, cockburn2007note, lehrenfeld2010hybrid}.

%% $H(\divergence)$-conforming finite elements for Navier-Stokes problems have been discussed and analyzed in \cite{cockburn2005locally, cockburn2007note}.

%% In \cite{cockburn2005locally} Cockburn, Kanschat and Sch\"otzau showed that finite element methods 
%% can only provide local conservation properties and energy-stability at the same time if numerical solutions are $H(\divergence)$-conforming.

Compared to Continuous Galerkin (CG) methods the number of degrees of freedom of Discontinuous Galerkin methods increases significantly. This drawdack is often outwayed by the advantages of the method. However, when it comes to solving linear systems Discontinuous Galerkin methods suffer most from drastically increased globally coupled degrees of freedom.
%from the additional degrees of freedom as the number of couplings, i.e. non-zero entries in a system matrix, increases even more drastically.
An approach to compensate for this is the concept of Hybridization where additional unknowns on element interfaces are introduced. This increases the number of degrees of freedom, but introduces two advantages. First, the global couplings are reduced and secondly, the structure of the couplings allows to apply static condensation for the element unknowns.
The concept has originally been introduced in the context of mixed finite element methods, cf. \cite{brezzi2012mixed}. In the last decade Hybridization has also been applied to Discontinuous Galerkin method for a variety of problems, cf. \cite{egger2008hybrid,nguyen2009implicit,cockburn2009unified,cockburn2011analysis,cesmelioglu2013analysis}. 
We also mention, that recently similar concepts are also known in the literature under the name ``Hybridized Weak Galerkin'' methods \cite{zhai2015hybridized} where a slightly different framework is used to derive the methods.
With the emphasis to treat general polyhedral meshes ``Hybrid High-Order'' (HHO) methods \cite{dipietro14,di2015review,cockburndipietroern2015} have recently been introduced where a combination of Hybrid (or hybridized) methods and higher order spaces is considered. 

Concerning HDG discretizations for incompressible flows, we also mention the papers \cite{egger2012hp,cockburn2011analysis,nguyen2010hybridizable}. Further, different approaches to implement exact incompressible finite element solutions to the Stokes problem using Hybridization have been investigated in \cite{carrero2006hybridized,cockburn2005incompressible,cockburn2005bincompressible}.
In \cite{stenberg2011} the Stokes-Brinkman problem, the Stokes problem with an additional zero order term, has been discretized using a hybridized $H(\divergence)$-conforming finite element formulation.  

An interesting and often praised aspect of Hybrid mixed finite element methods is the fact that a post-processing step can be used to reconstruct interior unknowns of an increased higher order accuracy for elliptic problems \cite{brezzi2012mixed}. The same is also possible for HDG methods in mixed formulation which is another advantage over conventional DG methods, see e.g. \cite{cockburn2009unified}. For elliptic problems an accuracy of order $k+2$ in the volume can be achieved if polynomials of order $k$ are used on the element interfaces. One main contribution in this paper is the introduction of a projection operator which achieves the same effect. In \cite{lehrenfeld2010hybrid} we already discussed this operator which has recently also been addressed under the name ``reduced stabilization'' in \cite{qiu2013hdg,oikawa2014hybridized,oikawa2015reduced,qiu2015superconvergent}. In these papers the construction of the corresponding operator is only possible in two dimensions using special integration rules. We explain how the operator can be implemented in a more general setting.

The efficiency and practicability of HDG methods is rarely addressed in the literature. Interesting exception are \cite{kirby2012cg,huerta2013efficiency} where the computational cost of the method is compared to CG and DG methods.

%% In contrast to most of the mentioned HDG formulations we aim for a formulation in terms of velocity and pressure ... no post-processing..proposed in \cite{lehrenfeld2010hybrid} and also pursuit in \cite{qiu2013hdg,oikawa2014hybridized,oikawa2015reduced}.

The nature of diffusive (viscous) and convective terms appearing in convection-diffusion-type problems have a substantially different character. Discretizations of problems involving only one of the two mechanisms would typically lead to different \emph{spatial and temporal} discretizations. It is therefore often desirable to consider operator-splitting time discretization schemes which allow for a separate treatment of the different operators as in \cite{ascher1995implicit, ascher1997implicit,kanevsky2007application,strang1968construction,maday1990operator, chrispell2007fractional}. In these papers the spatial discretization for the stiff (diffusion/Stokes) and the non-stiff (convection) operators is typically the same. We consider a different treatment of the operators with respect to their \emph{temporal and spatial} discretization. 
% Without going into the details of these operator splitting schemes, we consider such a method and refer to the literature for details.
% In \cite{ascher1995implicit, ascher1997implicit,kanevsky2007application} Implicit-Explicit (IMEX) schemes are discussed, in \cite{strang1968construction,maday1990operator} Operator-Integration-Factor Splitting schemes are presented and in \cite{glowinski2003finite} other time integration schemes like the fractional step method are considered.

\subsection{The concept: Efficiency through operator splitting}
%\todo[inline]{introduce notation for the operators. $M,A,C$?}
A crucial ingredient in the considered discretization is the fact that the convection and the Stokes problem are separated by means of an operator-splitting method. 
The operator-splitting method is chosen such that only operator evaluations of the convection operator are required so that the time integration scheme is explicit in terms of the convection operator. This is often affordable as the time step restrictions following from the convection operator are the least restrictive ones in the Navier-Stokes problem. Moreover, the convection operator is non-linear, s.t. the set up of linear systems of equations and corresponding preconditioners or solvers would have to take place every time step which renders implicit approaches for the convection very expensive. The Stokes operator is dealt with implicitly, i.e. it appears in linear systems in every time step. Due to the differential-algebraic structure of the Navier-Stokes equations this is necessary w.r.t. the pressure and the incompressibility constraint. As completely explicit handling of the viscosity terms would introduce severe time step restrictions it is also advisable to treat viscous forces implicit. Note that the Stokes operator is time-independent such that the setup of linear systems and preconditioners or solvers can be done once and re-used in every time step.

The spatial discretizations for the different operators are designed differently as the treatment of the convection term can be optimized for operator evaluations while the discretization of the Stokes operator has to provide an efficient handling of implicit solution steps, i.e. the solution of corresponding linear systems. This different treatment reflects in the use of two different finite element spaces which are used: An $H(\divergence)$-conforming Hybrid DG space for the velocity-pressure-pair for solving Stokes problems and a DG space for handling convection. 
%Between these spaces transfer operations are necessary. The spaces and the transfer operations that can easily be provided are presented in the next section. Later, in section \ref{sec:opsplitting}, we present operator-splitting methods, which can exploit the considered splitting of spatial operators and which only rely on the transfer operations provided.

\subsection{Main contributions and structure of this paper}
In this paper we introduce a new discretization method for the incompressible Navier-Stokes equations. The discretization is based on a decomposition of the problem into the (unsteady) Stokes problem and a hyperbolic transport problem. This decomposition reflects in the use of appropriate operator splitting methods in the time discretization and the use of different finite element spaces for the different spatial operators.

While we use a rather standard DG formulation for the spatial discretization of the hyperbolic transport problem, we consider a new HDG formulation for the spatial discretization of the Stokes-type problem. 
We use $H(\divergence)$-conforming functions of element-wise polynomials of degree $k$ for the velocity and discontinuous pressure functions of degree $k-1$. 
Due to the introduction of additional unknowns of degree $k-1$ on the facets for the approximation of the tangential trace of the velocity static condensation allows to reduce the unknowns and their couplings significantly. The resulting globally coupled unknowns correspond to the approximation of the normal trace of the velocity on the facets by order $k$ functions, the tangential trace of the velocity by order $k-1$ functions and the approximation of the pressure field by piecewise constants.

The main contributions of this paper are:
\begin{enumerate}
\item Introduction of a new $H(\divergence)$-conforming high order accurate HDG method with a \emph{projected jumps} formulation (also known as \emph{reduced stabilization} or \emph{reduced-order} HDG) for the solution of Stokes-type problems.
\item Presentation of a combined spatial discretization for the Navier-Stokes equations based on a standard Upwind DG formulation for the hyperbolic transport problem with the new HDG method for Stokes-Brinkman problems. 
\item Discussion of operator-splitting time integration schemes which restrict solutions of linear systems to Stokes-Brinkman problems. 
\end{enumerate}
In section \ref{sec:DGHDG} we discuss the discretization of spatial operators. The discussion is divided into two parts, the discretization of the Stokes-Brinkman problem, the problem which involves all relevant spatial operators except for the convection and the discretization of the convection operator. For the former part we consider an HDG discretization, for the latter part a standard DG discretization. 
%The HDG discretization allows a very efficient treatment of arising linear systems while the DG discretization is very suitable for an efficient treatment of operator evaluations. 
In section \ref{sec:opsplitting} operator-splitting time integration schemes are discussed which are tailored for such a situation. 
Finally, in section \ref{sec:numex} we give numerical examples which demonstrate the accuracy and the performance of the method.
% before we conclude the paper.

\subsection{Preliminaries}\label{hdgnse:notation}
Before describing the methods, we introduce some basic notation and assumptions:
$\Omega$ is an open bounded domain in $\mathbb{R}^d$ with a Lipschitz boundary $\Gamma$. It is decomposed into a shape regular partition $\mathcal{T}_h$ of $\Omega$ consisting of elements $T$ which are (curved) simplices, quadrilaterals, hexahedrals, prisms or pyramids. For ease of presentation we assume that all elements are not curved and consider only homogeneous Dirichlet boundary conditions. The element interfaces and element boundaries coinciding with the domain boundary are called \emph{facets}. The set of those facets $F$ is denoted by $\mathcal{F}_h$ and there holds $\bigcup_{T\in \mathcal{T}_h} \partial T = \bigcup_{F\in \mathcal{F}_h} F $. 
%The set of element interfaces, i.e. facets which are not on the boundary, is denoted by $\mathcal{F}^{int}_h$, whereas the set of boundary facets is $\mathcal{F}^{ext}_h$. 
%The inflow part of the domain or element boundaries are denoted by $\Gamma_{in}$ and $\partial T_{in}$, respectively.

In the sequel we distinguish functions with support only on facets indicated by a subscript $F$ and those with support also on the volume elements which is indicated by a subscript $T$. Compositions of both types are used for the HDG discretization of the velocity which is denoted by $u = (\HdivVar,u_F)$.

For vector-valued functions the superscripts $t$ denotes the application of the tangential projection: $v^t = v - (v\! \cdot\! n)\cdot n \in \mathbb{R}^d$. The index $k$ which describes the polynomial degree of the finite element approximation at many places through out the paper is an arbitrary but fixed positive integer number. 

%For the discretization of the spatial operators we mainly discuss linear, bilinear and trilinear forms. 
We identify finite element functions $u \in X_h$ with their representation in terms of coefficient vectors $\mathbf{u} \in \rr^{N_X}$, s.t. $u = \sum_{i=1}^{N_X} \mathbf{u}_i \varphi_i$ for a corresponding basis $\{\varphi_i\}$ of $X_h$ and $N_X = \mathrm{dim}(X_h)$. A (generic) bilinear form $\bilinearform{G}: X_h \times Y_h \rightarrow \rr$ is identified with the matrix $\mathbf{G} \in \rr^{N_Y \times N_X}$, s.t. $\mathbf{G}_{i,j} = \bilinearform{G}(\varphi_j^X,\varphi_i^Y)$.

\section{DG/HDG spatial discretization} \label{sec:DGHDG}
In this section we introduce the spatial discretization for the Stokes operator, the convection operator and transfer operations between both. 
First, we introduce the $H(\divergence)$-conforming Hybrid DG discretization of the Stokes-Brinkman problem, i.e. the stationary Navier-Stokes problem without convection and an additional zero order term in section in section \ref{hdgnse:stokes}. This discretization is improved significantly. Further on, a modification of the discretization using the idea of \emph{projected jumps} (also known as \emph{reduced stabilization} or \emph{reduced-order} HDG) is presented in section \ref{sec:modifications} including an explanation of how the projection operator can be realized. 
For the convection part of the Navier-Stokes problem we consider a DG discretization using standard approaches. This is discussed in section \ref{dg:convection}. As the discretization spaces for the Stokes part and the convection part are different, we present transfer operations between the spaces in section \ref{sec:transfer} which allow us to finally formulate the semi-discrete problem in section \ref{hdgnse:semidiscrete}.

\subsection{$H(\divergence)$-conforming HDG formulation of the Stokes-Brinkman problem} \label{hdgnse:stokes}
In this part we consider the discretization of the Stokes part of the Navier-Stokes problem. For simplicity we restrict the discussion to homogeneous Dirichlet boundary conditions. We present the method and elaborate on important properties. For an error analysis of the method we refer to \cite{lehrenfeld2010hybrid}.
The reaction term $\tau^{-1}$ corresponding to an inertia term stemming from an implicit time integration scheme is further incorporated. The resulting problem is known as the (stationary) Stokes-Brinkman problem:
\begin{equation}
\left\{
\begin{array}{r l l l}
\tau^{-1} u + \divergence( - \nu {\nabla} \VelVar  + p {I} ) &= {f} \quad &\mbox{ in } \Omega \\
   \divergence \ \! \VelVar &= 0 \quad &\mbox{ in } \Omega \\
%  				\VelVar &= 0 \quad &\mbox{ on } \Gamma_D \\
%   				(\nu {\nabla} \VelVar  - p {I}) \cdot {\normal} &= 0 \quad &\mbox{ on } \Gamma_{out} \\
\end{array} \right.
\label{eq:stokes}
\end{equation}
We first introduce the finite element spaces, followed by the definition of the bilinear forms corresponding to the involved operators.
%present the basic HDG discretization before we discuss an improvement of the method in section \ref{sec:modifications}.

%\subsection{Finite element spaces for the Stokes problem} \label{sec:spaces}
\subsubsection{$H(\divergence)$-conforming Finite Elements for Stokes.}\label{hdgnse:hdivfe}
Following \cite{cockburn2005locally} a DG formulation for the incompressible Navier Stokes equations which is locally conservative and energy-stable at the same time has to provide discrete solutions which are exactly divergence-free. This can be achieved with a suitable pair of finite element spaces. We consider the use of $H(\divergence)$-conforming Finite Element spaces for the velocity $\VelVar$. 
$H(\divergence)$-conformity requires that every discrete function $\VelVar_h$ is in 
\begin{equation*}
 H(\divergence,\Omega) = \{ \VelVar \in [L_2(\Omega)]^d:
\divergence\ \VelVar \in L_2(\Omega) \}
\end{equation*}
We use piecewise polynomials, so that on each element the functions are in $H(\divergence,T)$. For global conformity, continuity of the normal component is necessary, resulting in 
\begin{equation}
 \HdivSpace := \{ \HdivVar \in \prod_{T \in \mathcal{T}_h}  [\mathcal{P}^k(T)]^d, \ \jumpleft \HdivVar\! \cdot\! n \jumpright_F = 0 \ \forall \ F \in \mathcal{F}_h \} \subset H(\divergence,\Omega), \quad N_{\HdivSymb} := \mathrm{dim}(\HdivSpace),
\end{equation}
with $\jumpleft \cdot \jumpright_F$ the usual jump operator and $\mathcal{P}^k$ the space of polynomials up to degree $k$. We refer to \cite{brezzi2012mixed} for details on the construction of $H(\divergence)$-conforming finite elements such as the Brezzi-Douglas-Marini finite element. \\
The appropriate Finite Element space for the pressure is the space of piecewise polynomials which are discontinuous and of one degree less:
\begin{equation}
 \PressureSpace := \prod_{T \in \mathcal{T}_h} \mathcal{P}^{k-1}(T), \quad N_{\PressureSymb}  := \mathrm{dim}(\PressureSpace).
\end{equation}
This velocity-pressure pair $\HdivSpace$/$\PressureSpace$ fulfills 
\begin{equation}
\divergence ( \HdivSpace ) = \PressureSpace. 
\end{equation}
The crucial point of this choice of the velocity-pressure pair is the property:\\
If a velocity $\HdivVar \in \HdivSpace$ is weakly incompressible, it is also strongly incompressible:
\begin{equation} \label{eq:stronginc}
\int_{\Omega} \divergence( \HdivVar) \PressureTest \, dx = 0 \ \forall q_h \in \PressureSpace \quad \Leftrightarrow \quad \divergence(\HdivVar)=0 \text{ in } \Omega.
\end{equation}
% i.e. $\divergence(u_h)=0$ holds pointwise for $u_h \in \HdivSpace$.
% Additionally it satisfies an \emph{inf-sup} condition, called the \emph{LBB}-condition as we will see in the next section. \\
The benefit of \eqref{eq:stronginc} is twofold: First, it allows to show energy-stability for the incompressible Navier-Stokes equations, cf. section \ref{hdgnse:semidiscrete}. Secondly, error estimates for the velocity error can be derived which are independent of the pressure field, we comment on this in remark \ref{rem:linke} below.

As solutions of the incompressible Navier Stokes equations are $[H^1(\Omega)]^d \times L_2(\Omega)$-regular and tangential continuity is not imposed as an essential condition on the finite element space the discrete formulation has to incorporate the tangential continuity weakly. We do this with a corresponding (Hybrid) DG formulation for the tangential components across element interfaces 
%in section \ref{hdgnse:stokes}
%
\begin{remark}[Reduced spaces] \label{rem:red}
Due to \eqref{eq:stronginc} solutions of the discretized (Navier)-Stokes problem will be exactly divergence-free velocity fields. This a priori knowledge can be exploited. The basis for the space $\HdivSpace$ can be constructed in such a way such that we can discard certain higher order basis functions (with non-zero divergence) that will have no contribution. A corresponding basis is introduced in \cite{zaglmayr2006high} in the context of a set of higher order basis functions which fulfill an exact sequence property. The resulting space $\HdivSpace^{\text{red}}$ has $\divergence(\HdivSpace^{\text{red}}) = \PressureSpace^{\text{red}} := \{ \PressureVar|_{T} = const : T \in \mathcal{T}_h \}$ such that also most of the degrees of freedom of the pressure space can be discarded. The reduction of basis functions for $\HdivSpace$ and $\PressureSpace$ is explained in more detail in \cite[Chapter 2.2]{lehrenfeld2010hybrid}. 
%We note that the basis functions affected by the reduction are contained in the set of basis functions that can be locally eliminated by static condensation when it comes to solving linear systems involving the Stokes operator.
\end{remark}

\subsubsection{The HDG space for the velocity.} \label{sec:velspaces}
To (weakly) enforce continuity we apply a discontinuous Galerkin (DG) formulation such as the Interior Penalty method \cite{arnold1982interior}. However, to avoid the full coupling of degrees of freedom of neighboring elements, we introduce additional unknowns on the skeleton, the \emph{facet unknowns}, which represent an approximation of the tangential trace of the solution. 
The DG formulation is then replaced with a corresponding HDG formulation, s.t. degrees of freedom of neighboring elements couple only through the facet unknowns. We note that the resulting space is only ``hybrid'' in the sense of the tangential unknowns. Further, we do not consider a hybridization with respect to the pressure, cf. also remark \ref{rem:comparison}. 
\begin{figure}[h]
  \begin{center}
    \setlength{\unitlength}{3cm}
    \begin{tabular}{c c c}
      \multirow{3}{*}
      {
      \raisebox{\ht\strutbox-1.5\totalheight}{
      \begin{picture}(1.5,0.4)
        \includegraphics[width=3.5cm]{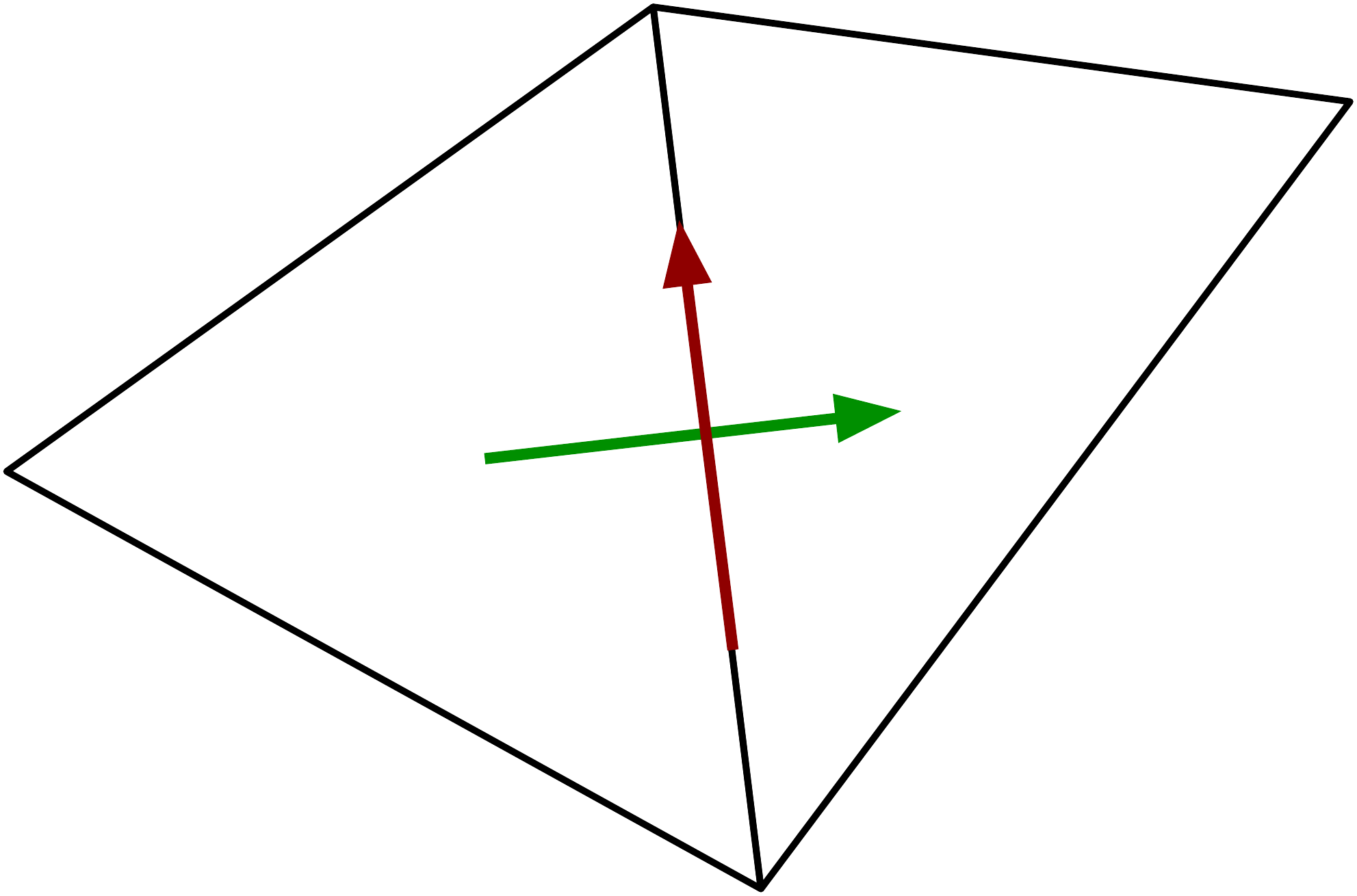}
      \end{picture}
      }
      }
      & \includegraphics[width=3.5cm]{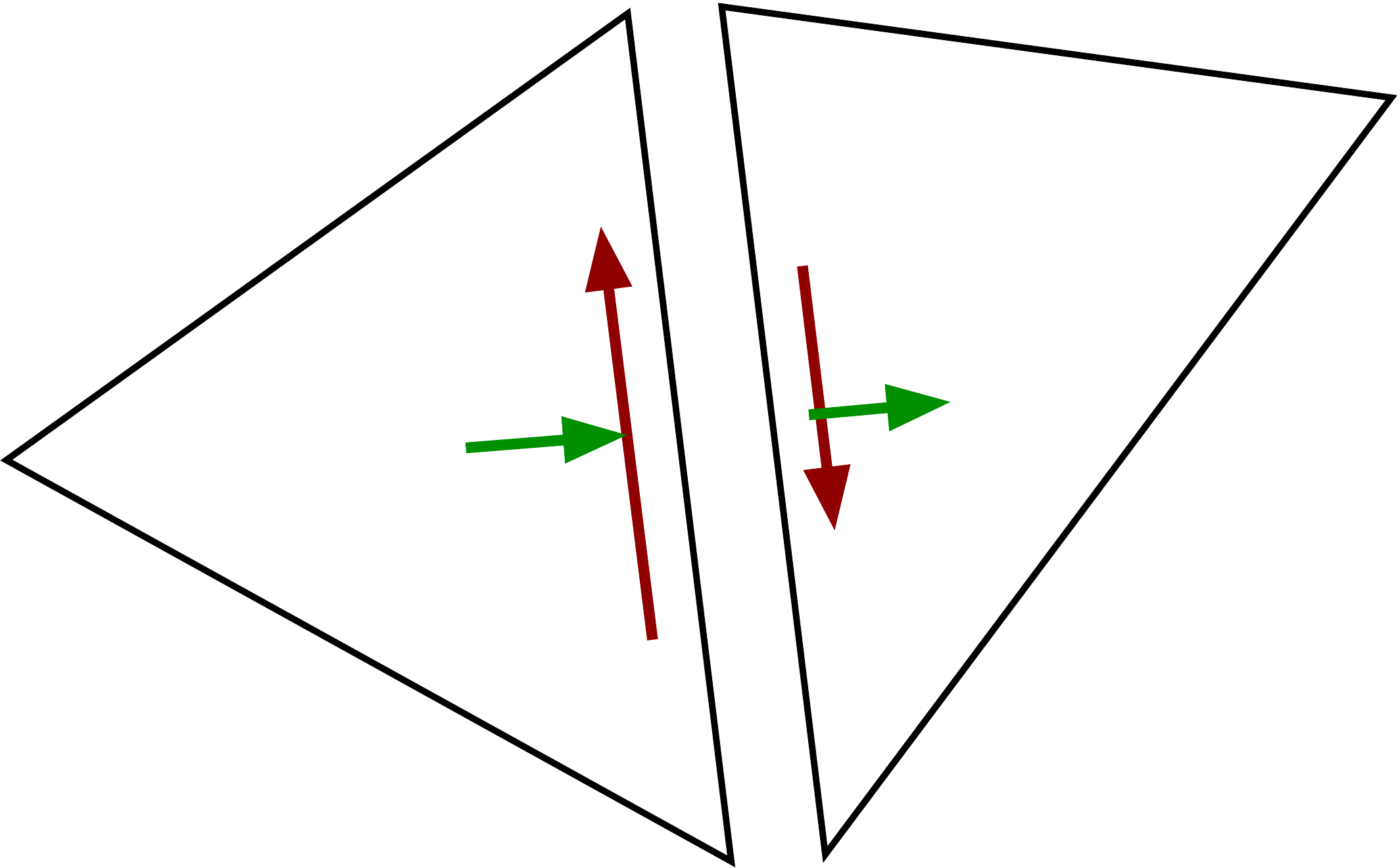}
      & \includegraphics[width=3.5cm]{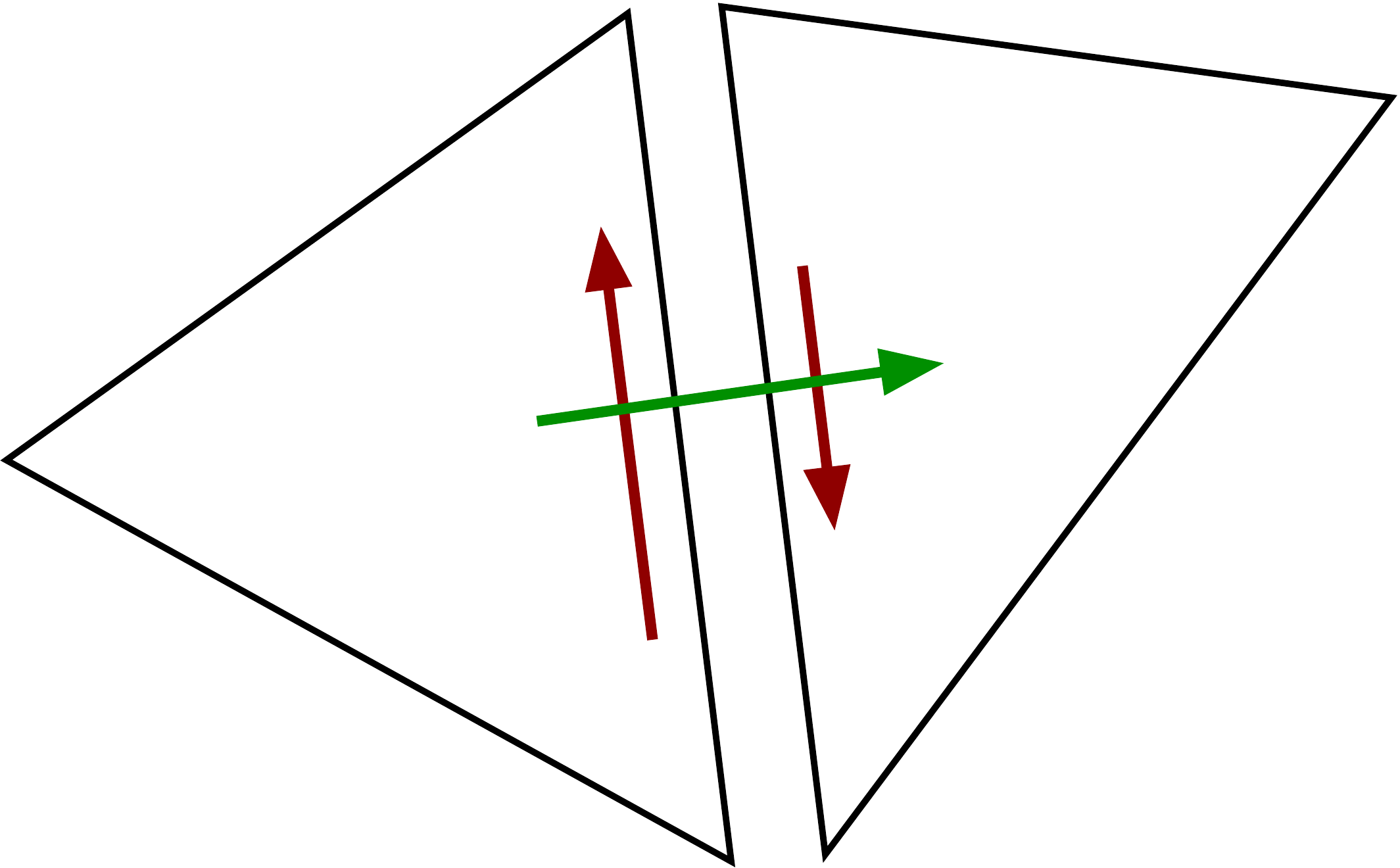}
      \\
      & DG 
      & $H(\divergence)$-conforming DG 
      \\
      & \includegraphics[width=3.5cm]{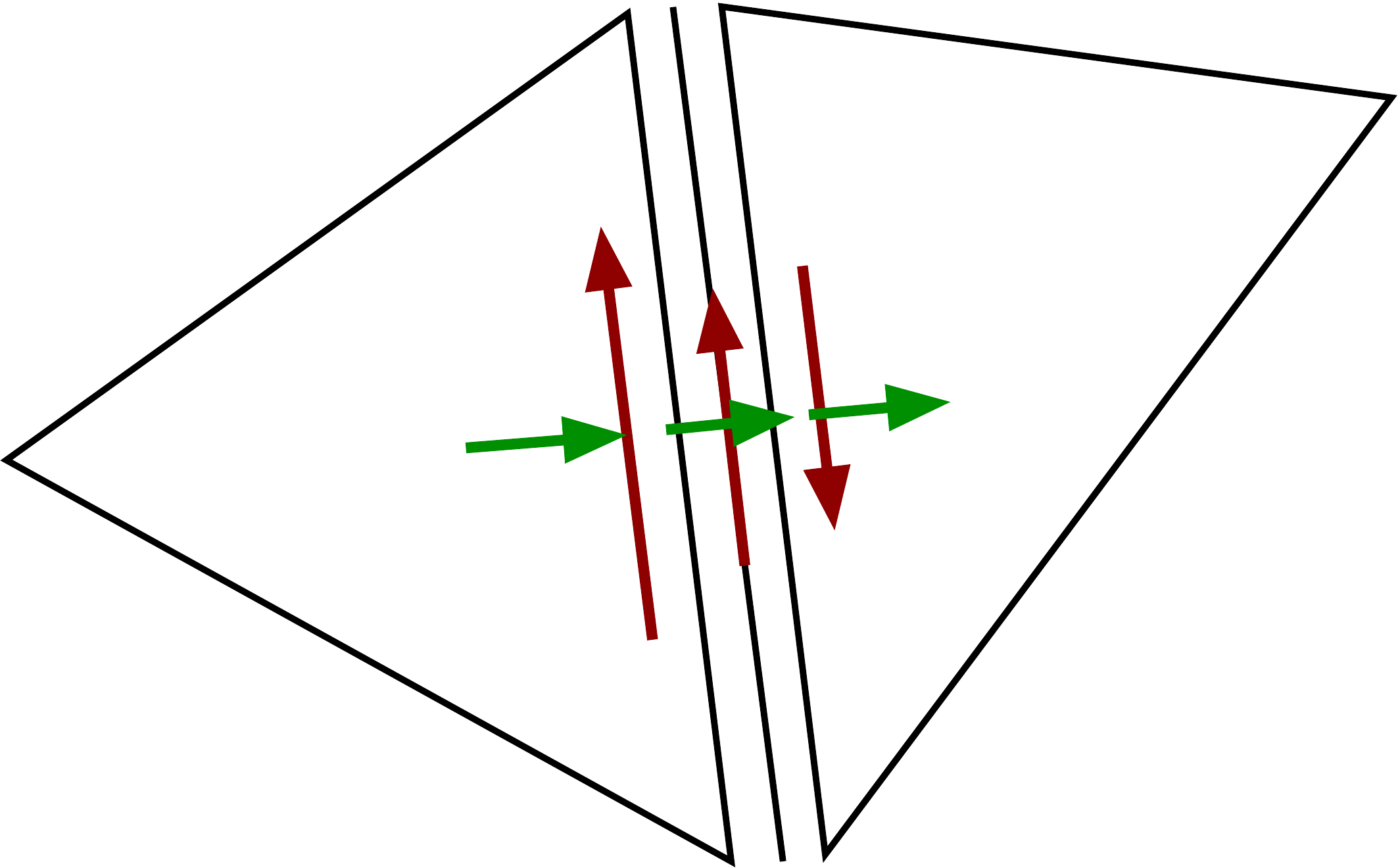}
      & \includegraphics[width=3.5cm]{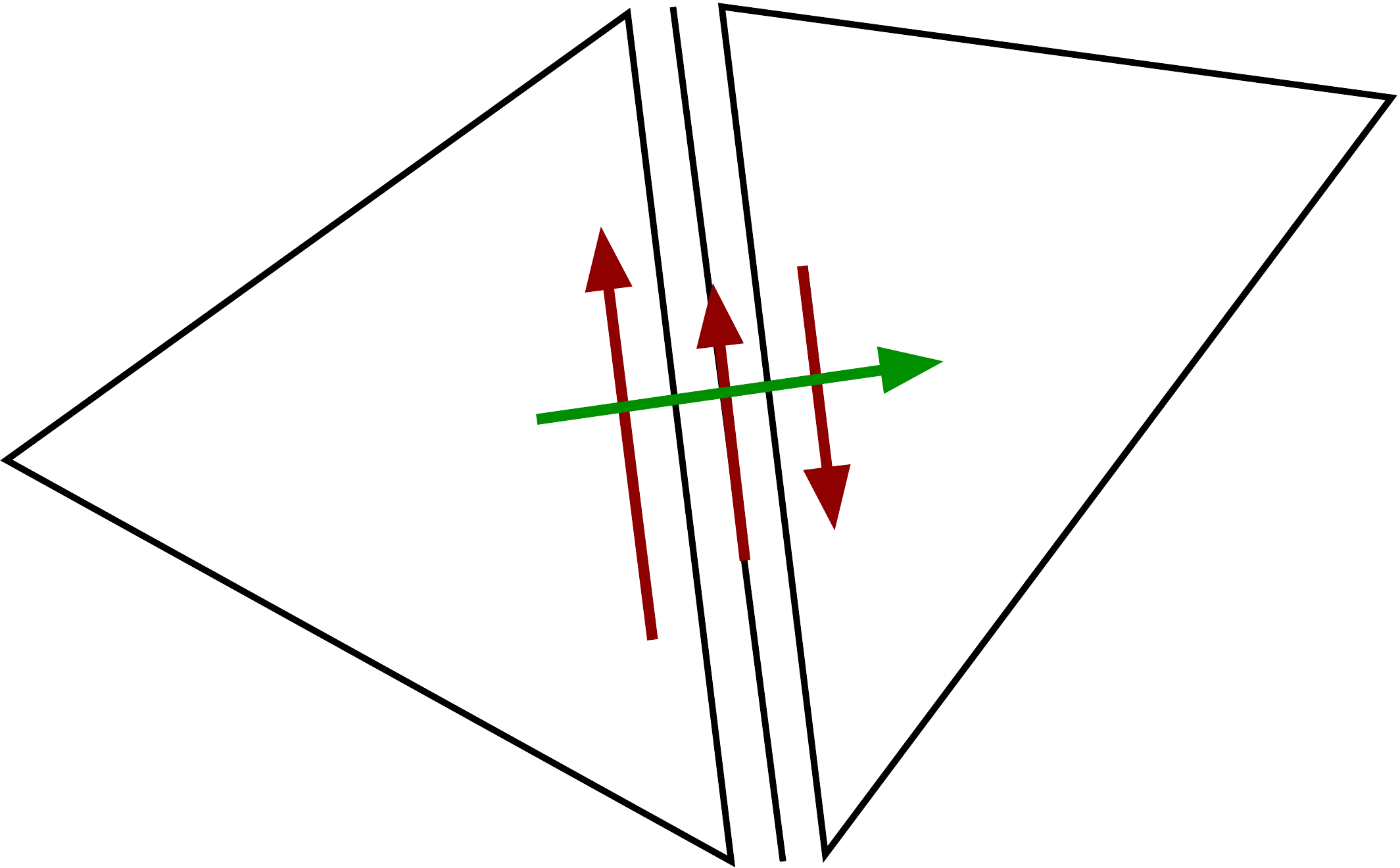}
      \\
      $H^1$-conforming (CG) 
      & HDG
      & $H(\divergence)$-conforming HDG 
    \end{tabular}
  \caption{Tangential and normal continuity for different finite element spaces for the velocity.}
  \label{fig:nseschemes}
  \end{center}
\vspace*{-0.5cm}
\end{figure}

\noindent 
As normal continuity is already implemented in $\HdivSpace$ we only need a DG enforcement of continuity in the tangential direction and hence only introduce the facet unknown for the tangential direction of the trace:    
% To this end, we introduce the space of facet unknowns:
\begin{equation}
\FacetSpace := \{ \FacetVar \in \prod_{F\in \mathcal{F}_h} [\mathcal{P}^k(F)]^d, \ \FacetVar \cdot  n = 0\ \}, \quad N_{\FacetSymb} := \mathrm{dim}(\FacetSpace).
\end{equation}
Functions in $\FacetSpace$ have normal component zero. For the discretization of the velocity field we use the composite space 
\begin{equation}
  \VelSpace := \HdivSpace \times \FacetSpace, \quad N_{\VelSymb} := \mathrm{dim}(\VelSpace).
\end{equation}

% \begin{remark}[$\VelVar$]
% \todo[inline]{At many places in the remainder we consider $\VelVar \in \VelSpace$ even if the facet contribution..}
% \end{remark}

\begin{remark}[The role of the facet space] \label{rem:rolefacet}
We note that the space $\FacetSpace$ is only introduced to allow for a more \emph{efficient handling of linear systems}. In section \ref{sec:numexlinsys} we consider a numerical example where a vector-valued Poisson problem is considered for the HDG space $\VelSpace$ in comparison to other ($H(\divergence)$-conforming) DG spaces to illustrate the impact of Hybridization. 
In case that only explicit applications of discrete operators are used, the introduction of the space $\FacetSpace$ entails no advantages. 
In this paper, facet variables appear only in the discretization of the viscous forces. Although the facet variable has no contribution in the discretization of inertia and pressure forces or the incompressibilty constraint we define the corresponding operations for $\VelVar \in \VelSpace$ instead of $\HdivVar \in \HdivSpace$ to simplify the presentation.
\end{remark}

\subsubsection{Viscous forces.} \label{sec:discvisc} 
For ease of presentation we consider a hybridized version of the Interior Penalty method for the discretization of the viscous term. In remark \ref{rem:ip} we comment on alternatives. 
In the hybridized version of the Interior Penalty method, the usual jump across element interfaces of the Interior Penalty method is replaced with jumps between element interior and facet unknown (in tangential direction) $\jumpleft u^t \jumpright = \HdivVar^t - u_F$ from both sides of a facet. 
%Note that due to normal-continuity $\jumpleft \HdivVar \cdot n \jumpright = 0$, $\HdivVar \in \HdivSpace$ the normal component of the jump is zero. 
The bilinear form corresponding to the HDG discretization of viscous forces is
\begin{equation} \label{eq:blfA}
\begin{array}{rl}
\displaystyle \bilinearform{A}(\VelVar,{v}) := & \displaystyle \sumoverallelements \int_{T} \nu {\nabla} {\HdivVar} \! : \! {\nabla} {\HdivTest} \ d {x} - \int_{\partial T} \nu \frac{\partial {\HdivVar}}{\partial {\normal} }  \jumpleft {v}^t \jumpright \ d {s} \\
& \displaystyle - \int_{\partial T} \nu \frac{\partial {\HdivTest}}{\partial {\normal} }  \jumpleft \VelVar^t \jumpright \ d {s}
 + \int_{\partial T} \nu \frac{\stab}{h}  \jumpleft \VelVar^t \jumpright  \jumpleft {v}^t \jumpright \ d {s}, \quad  u\!=\!(u_T,\!u_F), v\!=\!(v_T,\!v_F) \!\in\! U_h \hspace*{-0.5cm}
\end{array}
\end{equation}
In this bilinear form the four terms have different functions. While the first two terms ensure consistency (in the sense of Galerkin orthogonality) the third and fourth term are tailored to ensure symmetry (adjoint consistency) and stability, respectively. Due to continuity of the solution ($\jumpleft \cdot^t \jumpright = 0$) the latter terms also preserve consistency. For a more detailed introduction we refer to \cite[section 2.3.1]{lehrenfeld2010hybrid}.
With respect to the discrete norm
\begin{equation*}
  \brokenHnormIIleft \VelVar \brokenHnormIIright^2 := \sumoverallelements
\left\{  \Vert {\nabla} {\HdivVar} \Vert_T^2  + \frac{1}{h} \Vert \jumpleft
\VelVar^t  \jumpright \Vert_{\partial T}^2  + h \left\Vert \frac{\partial u_T}{\partial n} \right\Vert_{\partial T}^2  \right\}
\end{equation*}
which is an appropriately modified version of the discrete norm typically used in the analysis of Standard Interior Penalty methods, the bilinearform $\bilinearform{A}$ is (for a sufficiently large $\stab$) consistent, bounded and coercive.
%, cf. \cite{lehrenfeld2010hybrid}.

%
%
%If we compare this with a standard approach to discretize the vector-valued laplacian with the use of discontinuous Finite Element spaces, namely the Interior Penalty method (see e.g. \cite{unifiedDG}) you get the same formula except that the adjacent element's tangential value take the role of $\HdivTest$ and $\HdivVar$ within the jump operator. \\
%We briefly discuss the structure induced by the HDG formulation in contrast to DG formulations. 
The coupling through facet unknowns enforces the same kind of continuity as in the Interior Penalty DG formulation while preserving the following structure in the sparsity pattern: element unknowns are only coupled with unknowns associated with the same element or unknowns associated with aligned facets.
% This improvement can have a dramatic effect on the complexity of linear systems, cf. remark \ref{rem:rolefacet} and the numerical example in section \ref{sec:numexlinsys}.

\begin{remark}[Interior Penalty and alternatives] \label{rem:ip}
  A drawback of the (Hybrid) Interior Penalty method is the fact that the stabilization parameter $\stab$ depends on the shape regularity of the mesh. Often $\stab$ is chosen on the safe side, but as was pointed out in \cite{schoberl2013domain} the condition number of arising linear systems increases with $\stab$. %, s.t. you don't wont to choose $\stab$ to large. 
In the same paper a hybridized variant of the Bassi-Rebay stabilization method (cf. \cite{bassi1997high,bassirebay97b,brezzietal99}) for a scalar Poisson equation has been proposed, see also remark \ref{rem:precond}. Such a variant is used in the numerical examples, but as it does not have any further consequences for the remainder of this work, we stick to the well-known (Hybrid) Interior Penalty method for ease of presentation.
\end{remark}

\subsubsection{Mass bilinear form.} \label{sec:mass}
The HDG mass matrix is defined as
\begin{equation}
  \bilinearform{M}^{\VelSymb}(\VelVar, \VelTest) := \int_{\Omega} \HdivVar \HdivTest \, dx, \qquad \qquad\VelVar=(\HdivVar,\FacetVar), \VelTest=(\HdivTest,\FacetTest) \in \VelSpace.
\end{equation}
Note that we defined the mass matrix for $\VelVar \in \VelSpace$ although $\FacetVar$ has no contribution, cf. remark \ref{rem:rolefacet}.

\subsubsection{Pressure force and incompressibility constraint.} \label{sec:pressure}

For the pressure force and the incompressibility constraint we define the bilinearform
\begin{equation}
\bilinearform{D}({\VelVar},\PressureVar) := \sumoverallelements - \int_{T} \PressureVar \
\divergence \ \! {\HdivVar} \ d {x} \quad \text{ for } \VelVar=(\HdivVar, \FacetVar) \in \VelSpace, \ \PressureVar \in \PressureSpace
\end{equation}
for which the \emph{LBB-condition}
\begin{equation}\label{eq:lbb}
  \sup_{\VelVar \in \VelSpace} \frac{\bilinearform{D}(\VelVar,\PressureVar)}{\brokenHnormIIleft \VelVar \brokenHnormIIright} \geq c_{{}_{LBB}} \Vert \PressureVar \Vert_{L^2}, \quad \forall \ \PressureVar \in \PressureSpace
\end{equation}
holds true for a $c_{{}_{LBB}}$ independent of the mesh size $h$, cf. \cite[Proposition 2.3.5]{lehrenfeld2010hybrid}.

\begin{remark}[Robustness in polynomial degree $k$] \label{rem:probust}
In numerical experiments we observed that the LBB constant $c_{{}_{LBB}}$ in \eqref{eq:lbb} is also robust in $k$. To this end we computed the LBB constant of one element with Dirichlet boundary conditions and the condition $\int_{\Omega} p \, dx=0$ on the pressure. The results are shown in Table \ref{tab:lbb} for varying the polynomial degree $k$. 
\begin{table}
\begin{center}
\begin{tabular}{rcccc}
\toprule
 $k$ 
& 4  & 8  & 16  & 32
\\
\midrule
quadrilateral
& $0.305$ 
& $0.313$ 
& $0.315$ 
& $0.315$ 
% & $0.093$ 
% & $0.098$ 
% & $0.099$ 
% & $0.099$ 
\\
triangle
& $0.167$ 
& $0.190$ 
& $0.201$ 
& $0.205$ 
% & $0.028$ 
% & $0.036$ 
% & $0.040$ 
% & $0.042$ 
\\
% tetrahedron
% & $0.011$ 
% & $0.011$ 
% & $0.0xx$ 
% & $0.0xx$ 
% \\
\bottomrule
\end{tabular} 
\end{center}
\vspace*{-0.4cm}
\caption{LBB constant in dependence of $k$ for a single element.}\label{tab:lbb}
\vspace*{-0.4cm}
\end{table}
From the boundedness of the LBB constant on one element and the $h$-robustness in \eqref{eq:lbb}, the robustness in $h$ and $k$ on the global spaces follows. This will is in the forthcoming master's thesis of Philip Lederer.
\end{remark}

\subsubsection{Discretization of the Brinkman-Stokes problem} \label{sec:hdgdisc}
With the introduced discretizations of the spatial operators the discrete Brinkman-Stokes problem can be written as:\\
Find $\VelVar \in \VelSpace$ and $p \in \PressureSpace$, s.t.
\begin{eqnarray}\label{eq:discbrinkman}
\displaystyle
\left\{
\begin{array}{crcrl}
\tau^{-1} \bilinearform{M}(\VelVar,{v}) + \bilinearform{A}(\VelVar,{v}) &+ \ \ \bilinearform{D}({v},p) 	& = &
\langle {f},v \rangle 	& \forall \ v \in \VelSpace,\\
 & \bilinearform{D}(\VelVar,q)	& = & 0		&
\forall \ q \in \PressureSpace.
\end{array}
\right.
\end{eqnarray}
Due to coercivity of $\bilinearform{A}$ (respectively $\tau^{-1}\bilinearform{M}+\bilinearform{A}$), the LBB-condition of $\bilinearform{D}$ and consistency and continuity of all bilinear forms Brezzi's famous theorem (see \cite{brezzi2012mixed}) can be applied to obtain optimal order a priori error estimates. For a discussion of the coupling structure of this discretization we refer to Remark \ref{rem:comparison} below.
%%% TODO: put estimate here or not???
\begin{remark}[Pressure-independence of velocity error] \label{rem:linke}
Classical error estimates for mixed problems result in error estimates which are formulated in the compound norm of the velocity and pressure space. This has the disadvantage that the discretization error in the velocity depends on the approximation error in the pressure. As was pointed out in \cite{linke2014role}, ideally this should not be the case. Due to the fact that discrete solutions to \eqref{eq:discbrinkman} are exactly divergence-free an error estimate for the velocity field which does not involve the pressure can be derived, cf. \cite[Lemma 2.3.13]{lehrenfeld2010hybrid}.
\end{remark}

\subsection{Projected jumps: An enhancement of the HDG Stokes discretization}\label{sec:modifications}

The proposed HDG formulation for viscous forces in section \ref{sec:discvisc} can also be derived as a hybrid mixed method with a modified flux. This is done in detail in \cite{cockburn2009unified} (see also \cite[section 1.2.2]{lehrenfeld2010hybrid} or \cite{nguyen2010hybridizable}). In that setting the unknowns for the primal variable ($u_T$) are approximated with the same polynomial degree $k$ as the facet unknowns ($u_F$). Afterwards, in a postprocessing step approximations for the primal unknown, of one degree higher, $k+1$, are reconstructed in an element-by-element fashion. This approach ends up with a higher order approximation than the previously introduced HDG method considering the use of the same polynomial degree $k$ for the facet unknowns. This sub-optimality can be overcome by means of a projection operator which leads to the \emph{projected jumps} formulation. We first introduce the method and explain how the method can be implemented afterwards.

\subsubsection{Projected jumps: The method.}\label{sec:modifications1}
% Although the proposed formulation is optimally convergent, consistent and stable, there is, from a computational point of view, still potential for optimization. We consider 
% Two modifications are proposed in the next two sections. \\
% The first considers the construction of the basis functions where the (exact) divergencefree-constraint can be made use of to reduce the number of unknowns. The second approach reduces the number of facet unknowns by a slide modification of the discrete formulation. 
% \subsection{A ``projected jumps`` HDG formulation}\label{modifications:projectedjumps}
% \todo[inline]{blame postprocessing}
% How this disadvantage can be eliminated will be discussed in this section. 
%This approach is somewhat more direct than the ones we mentioned above as it does not involve any postprocessing step. 
%At the element boundaries, we replace jumps with $L_2$-projections of those jumps to polynomials of a degree lower. Thus the facet unknowns only have to be of order $k-1$, whereas the volume unknowns are of order $k$. 
%If this does not destroy the stability and consistency results of the original method, as it is the case here, it yields the same optimal order of convergence of the velocity with respect to the polynomial degree of the facet unknowns as the methods mentioned above. \\
The idea of the \emph{projected jumps} formulation is to reduce the polynomial degree  of the facet unknowns in $\FacetSpace$ to $k-1$,
\begin{equation}
\FacetSpace \rightarrow \FacetSpaceLO := \{ \FacetVar \in \prod_{F \in \mathcal{F}_h} [\mathcal{P}^{k-1}(F)]^d, \FacetVar \cdot  n = 0,\ \forall\ F \in \mathcal{F}_h\},
\end{equation}
while keeping the polynomial degree $k$ in $\HdivSpace$. In order to do this in a consistent fashion we have to modify the bilinearform $\bilinearform{A}$. 
First, we introduce the $L_2$ projection ${\proj}$ for a fixed facet $F \in \mathcal{F}_h$:
\begin{equation}
\displaystyle {\proj}: [\mathcal{P}^k(F)]^d \rightarrow [\mathcal{P}^{k-1}(F)]^d, \quad
\displaystyle  \int_F \left( {\proj} \ \! {f} \right) {v} \ d {x} = \int_F {f} \ \! {v} \ d {x} \quad \forall \ {v} \in [\mathcal{P}^{k-1}(F)]^d
\end{equation}
Due to the fact that 
$\frac{\partial{\HdivVar}}{\partial{{\normal}}} \in \mathcal{P}^{k-1}(F)$ on affine linearly mapped elements, the test functions ${\HdivTest}$
and ${\FacetTest}$ in the first integral on the element boundary in \eqref{eq:blfA} 
can be replaced by their $L_2(F)$ projections ${\proj} {\HdivTest}$
and ${\proj}\FacetTest$. If we additionally use
${\proj}\jumpleft \VelVar^t\jumpright$ instead 
of $\jumpleft \VelVar^t\jumpright$ to symmetrize and stabilize the
formulation in \eqref{eq:blfA} we end up with a modified version of the
Hybrid Interior Penalty method:\\
\begin{equation} \label{eq:blfAr}
\begin{array}{rl}
\displaystyle \bilinearform{A}^r(\VelVar,\VelTest) := & \displaystyle \sumoverallelements \int_{T} \nu {\nabla} {\HdivVar} \! : \! {\nabla} {\HdivTest} \ d {x} - \int_{\partial T} \nu \frac{\partial {\HdivVar}}{\partial {\normal} }  {\proj} \jumpleft {v}^t \jumpright \ d {s} \\
& \displaystyle- \int_{\partial T} \nu \frac{\partial {\HdivTest}}{\partial {\normal} } {\proj} \jumpleft \VelVar^t \jumpright \ d {s}
 + \int_{\partial T} \nu \frac{\stab}{h} {\proj} \jumpleft \VelVar^t \jumpright  {\proj} \jumpleft {v}^t \jumpright \ d {s}, \quad \VelVar, \VelTest \in \HdivSpace \times \FacetSpaceLO
\end{array}
\end{equation}
Note that the first two boundary integrals are just reformulated while only the last integral is really changed.
This modification preserves the important properties of $\bilinearform{A}$: It is still consistent and bounded. Further, a deeper look into the coercivity proof (see \cite{lehrenfeld2010hybrid}) reveals that coercivity can easily be shown in the modified (weaker) norm 
\begin{equation}
  \brokenHnormleft \VelVar \brokenHnormright^2 := \sumoverallelements
  \left\{  \Vert {\nabla} {\HdivVar} \Vert_T^2  + \frac{1}{h} \Vert \jumpleft {\proj} \ \!
    \VelVar^t \jumpright \Vert_{\partial T}^2 + h \left\Vert \frac{\partial \HdivVar}{\partial \normal} \right\Vert_{\partial T}^2  \right\}.
\end{equation}
Notice that, as we do not modify $\HdivSpace$, the normal component of $\HdivVar \in \HdivSpace$ is still a polynomial of degree $k$.

The idea of applying such a \emph{reduced stabilization} is also discussed and analyzed in \cite{oikawa2014hybridized,oikawa2015reduced}. However, only for the two-dimensional case the realization of the projection $\proj$ is discussed (by means of Gauss quadrature, cf. \cite[section 3.4]{oikawa2014hybridized}). 
In the next section a simple way to implement this projection operator is presented.

\begin{remark}[Interplay with operator-splitting]
If a Hybrid DG formulation for a problem involving diffusion (viscosity) and convection is applied, such a reduction of the facet unknowns is only appropriate if diffusion is dominating. 
One premise of the operator-splitting in this paper is that convection is not involved in implicit solution steps. Hence, the facet variable will never be involved in the discretization of convection, s.t.
even if the physical problem is convection dominated, the \emph{projected jumps} modification can be applied without loss of accuracy.
\end{remark}

\begin{remark}[Comparison to a fully hybridized formulation] \label{rem:comparison}
We note that the formulation in \eqref{eq:discbrinkman} (and the corresponding reduced formulation with $\mathcal{A}_h^r(\cdot,\cdot)$) is not a ``hybridized'' formulation in the usual sense, see e.g. \cite{cockburn2009unified}. Only for the tangential component of the velocity a new unknown is introduced, such that we only have a ``tangentially hybridized'' velocity discretization. The pressure is not hybridized. 
We comment on the resulting coupling of this HDG method and compare it to a HDG formulation where only Hybrid unknowns on the facet exist, for the velocity and the pressure. 
%The method in this paper is only ``hybrid'' in the tangential unknowns for the velocity which however suffices to obtain a global system which only involves velocity unknowns on the facets. 
The velocity unknowns in $U_h = W_h \times F_h$ (or $W_h \times \overline{F}_h$ in the reduced case) can be reduced to unknowns on the skeleton by static condensation.
The remaining \emph{facet unknowns} are the usual unknowns for the normal component of the $H(\divergence)$-conforming finite element space and the (tangential) unknowns in $F_h$ (or $\overline{F}_h$). 
Due to the pressure unknowns, static condensation of the formulation in \eqref{eq:discbrinkman} does not yield a global system which only involves facet unknowns, but also element unknowns for the pressure. 
However, except for the constant pressure all pressure unknowns can be eliminated by static condensation. We summarize: The globally reduced system only involves velocity unknowns on the facets and a constant pressure on each element.
Alternatively one could formulate a related HDG formulation which only involves facet unknowns for the velocity and the pressure. In such a formulation the pressure plays the role of the lagrange multiplier for the normal continuity and only for the (weak) tangential continuity velocity unknowns have to be added on the facet. To obtain an accurate (and exactly divergencefree) order $k$ approximation of the velocity and an order $k-1$ approximation of the pressure, the facet unknowns would be an order $k$ approximation of the pressure and an order $k$ (or order $k-1$ if projected jumps or corresponding postprocessing techniques are applied) approximation of the tangential velocities. Overall the number of unknowns would be smaller compared to the formulation proposed before by one (constant) pressure unknown per element. 
However, the price for this is the fact, that the Lagrange multiplier space (the pressure unknowns on the facets) is much larger and hence the structure of the arising linear system is different. 
It very much depends on the linear solver strategies which of the two Hybrid DG formulations gives the better perfomance. 
\end{remark}

\subsubsection{Projected jumps: Realization.}\label{implementations:projectedjumps}
The following way to implement \emph{projected jumps} needed for $\bilinearform{A}^r$ in \eqref{eq:blfAr} relies on an $L_2$-orthogonal basis for the facet functions in $\FacetSpace$. 
To obtain an $L_2$-orthogonal basis, we take a local coordinate system on each facet (spanned by $d-1$ aligned edges) and an $L^2$-orthogonal Dubiner basis (see \cite{dubiner1991spectral}) for each vector component. We  consider the three dimensional case, and denote the vectors of the local coordinate system by $e_1$ and $e_2$. Translation to the two dimensional case is then obvious. 

On each facet $F \in \mathcal{F}_h$, the space of (vector-valued) polynomials up to degree $k$, $\mathrm{span}(e_1,e_2) \cdot \mathcal{P}^k(F)$, can be split into the orthogonal subspaces
\begin{equation} \label{eq:projsplit}
 \mathcal{V}_{k-1} := \mbox{span}\{e_1,e_2\} \cdot \mathcal{P}^{k-1}(F) \quad \text{ and } \quad \mathcal{V}_{k-1}^\perp := \mbox{span}\{e_1,e_2\} \cdot [\mathcal{P}^{k}(F) \cap \left( \mathcal{P}^{k-1}(F) \right)^{\perp}].
\end{equation}
On each facet $F$ the facet unknown $\FacetVar$ can thus be written as $\FacetVar = \FacetVarLO + \FacetVarHO$ with unique $\FacetVarLO \in \mathcal{V}_{k-1}$, $\FacetVarHO \in \mathcal{V}_{k-1}^\perp$. We now replace the highest order function $\FacetVarHO$ with two copies $\FacetVarT$ and $\FacetVarTT$, each of these functions is associated with one of the two neighbouring elements $T$ and $T'$, the functions $\FacetVarT$ are only defined element-local. 
%We will use these functions to construct the projection method in an element-by-element fashion. 
$\FacetVarT$ can thus be eliminated after the computation of the element matrix and finally, only $\HdivVar \in \HdivSpace$ and $\FacetVarLO \in \FacetSpaceLO$ appear in the global system. 

$\FacetVarT$ is only responsible for implicitly realising the projection operator. 
To see this we 
%Next, we show that this results in an implementation of the \emph{projected jump}, 
now consider one facet $F \in \mathcal{F}_h$ and one of the neighboring elements $T \in \mathcal{T}_h$. We use the decomposition corresponding to \eqref{eq:projsplit} for trial and test functions
\begin{equation}
\begin{array}{l}
  \FacetVar = \FacetVarLO + \FacetVarT \mbox{ and } 
  \FacetTest = \FacetTestLO + \FacetTestT,
\end{array}
\end{equation}
with $\FacetVarLO, \FacetTestLO \in \mathcal{V}_{k-1}$ and $\FacetVarT, \FacetTestT \in \mathcal{V}_{k-1}^\perp$ where $\FacetVarT, \FacetTestT$ are only supported on $T$.
In \eqref{eq:discbrinkman} we choose the test function $v$ such that ${\HdivTest} = \FacetTestLO = 0$, $\FacetTestT \in \mathcal{V}_{k-1}^\perp$ on $T$ and $\FacetTestTT = 0$ on every other element $T'$. This yields
\begin{equation}
\int_F \nu \frac{\partial {\HdivVar}}{\partial {\normal}} {\FacetTestT} \ d {s} + \int_F \nu \frac{\stab}{h} (\HdivVar^t - \FacetVarLO - \FacetVarT) {\FacetTestT} \ d {s} = \int_F \nu \frac{\stab}{h} (\HdivVar^t - \FacetVarT) {\FacetTestT} \ d {s} = 0, \ \forall \ \FacetTestT \in \mathcal{V}_{k-1}^\perp.
\end{equation}
Hence, there holds ${\FacetVarT} = (I-{\proj})\HdivVar^t$ and we have
%This equation determines $\FacetVarT$ as 
\begin{equation}
% \ \ \Longrightarrow \ \
% \end{equation}
% Further we get: 
% \begin{equation}
\jumpleft \VelVar^t \jumpright = \HdivVar^t - \FacetVarLO - \FacetVarT = {\proj} (\HdivVar^t) - \FacetVarLO = {\proj}(\HdivVar^t - \FacetVarLO) = {\proj}\jumpleft \VelVar^t \jumpright
\end{equation}
We conclude that it is sufficient to consider the HDG bilinear form $\bilinearform{A}$ as before, where the local element matrices are computed according to $\FacetSpace$. For each element matrix the degrees of freedom corresponding to $\mathcal{V}_{k-1}^\perp$ are then eliminated forming a corresponding Schur complement. This yiels a final stiffness matrix which is only set up with respect to space unknowns of $\HdivSpace \times \FacetSpaceLO$.
%only has only the degrees of freedom of $\FacetSpaceLO$. 

\subsection{DG formulation for the convection}\label{dg:convection}

% \subsubsection{DG space for convection.} \label{dg:space}

For the discretization of the convection part we consider a standard DG finite element space:
\begin{equation}
\DGSpace := \{ u : u \in [\mathcal{P}^k(T)]^d \ \forall\, T \in \mathcal{T}_h \}, \quad N_{\DGSymb} := \mathrm{dim}(\DGSpace).
\end{equation}
We define the mass bilinear form
\begin{equation}
\bilinearform{M}^{\DGSymb}(\DGVar, \DGTest) :=  \int_{\Omega} {\DGVar} \ \!
{\DGTest} \ d{x}, \qquad \DGVar, \DGTest \in \DGSpace.
\end{equation}
Using an $L^2$-orthogonal basis on each element renders the associated mass matrix $\operator{M}^{\DGSymb}: \rr^{N_{\DGSymb}} \rightarrow \rr^{N_{\DGSymb}}$ diagonal. 
A stable spatial discretization of \eqref{eq:navierstokes} with respect to the convection part is achieved using a Standard Upwind DG trilinearform. We assume that the convection velocity $\HdivVar$ is  exactly divergence-free.
\begin{equation}
\begin{array}{rl}
    \displaystyle \bilinearform{C}(\HdivVar;\DGVar,\DGTest) := 
  \! \sum_T - \! \int_{T} \!\! \DGVar \otimes \HdivVar \!:\! \nabla \DGTest \ d{x} +
  \int_{\partial T} \!\! \HdivVar \! \cdot \! \normal \ \hat{\DGVar} \ \DGTest \ d {s}, \ \
  \HdivVar \in \HdivSpace, \ \DGVar, \DGTest \in \DGSpace.
%   \mbox{for } u,v,w \in \DGSpace \mbox{ or }
\end{array}
\end{equation}
where $\hat{\DGVar}$ denotes the upwind value $ \hat{\DGVar} = \lim_{\varepsilon \searrow 0} {\DGVar}({x}-\varepsilon \HdivVar({x})) $.
%
%We don't consider hybrid versions of the convection trilinearform as we will use time integration methods only which
%treat this part explicitely. Nevertheless if systems of equations involving $C_h$ are considered an hybridized version like it
%was considered in \cite{mixedhybriddg convdiff} might be of interest (see also \cite{diploma}). \\
Standard Upwind DG formulations are stable in the sense that with $\partial \Omega_{in} := \{ x\in \partial \Omega, \HdivVar(x) \cdot n < 0\}$ there holds
\begin{equation} \label{convstab}
\frac12 \int_{\partial \Omega_{in}} \!\!\!\!\! |\HdivVar\cdot n|\ w^2 \, ds + \bilinearform{C}({\HdivVar};\DGVar,\DGVar) \geq 0, \qquad \forall \ \DGVar \in \DGSpace, \ \VelVar \in \VelSpace, \ \divergence(\HdivVar) = 0.
\end{equation}
%Note that with $\VelSpace \subset \DGSpace$ the bilinear form $\bilinearform{C}(\HdivVar;\cdot, \cdot)$ can easily be extended to the HDG space $\VelSpace$.
\subsection{Transfer operations and embeddings}\label{sec:transfer}
%In this section we assume that no elements in $\mathcal{T}_h$ are curved, s.t. $\HdivSpace \subset \DGSpace$ and remark on the general case in remark \ref{rem:embcurv}.  
The discrete convection and Stokes operators have been defined on different spaces. 
%The spaces that the discretizations of the Stokes and convection operator are operating on are different. 
In order to combine both we introduce transfer operations between the (finite element) spaces to make both discretizations compatible. 
We restrict ourselves to two types of transfer operations which are based on embeddings:
\begin{itemize}
\item[$\transferop$:] With $\HdivSpace \subset \DGSpace$ we have
a canonical embedding of $\VelSpace$ in $\DGSpace$ with the embedding operator $\mathcal{I}: (\HdivVar,\FacetVar) \in \VelSpace \rightarrow \HdivVar \in \DGSpace$. Note that the corresponding operation in terms of coefficient vectors, denoted by $\transferop: \rr^{N_{\VelSymb}} \rightarrow \rr^{N_{\DGSymb}} $, is not an identity.
\item[$\transferopI$:] The embedding operator $\mathcal{I}$ implies the canonical embedding $\mathcal{I}': \DGSpace' \rightarrow \VelSpace'$, which maps functionals 
$ \mathcal{I}': f \in \DGSpace' \rightarrow [ (\HdivVar,\FacetVar) \in \VelSpace \rightarrow f(\HdivVar) ] \in \VelSpace'$. The corresponding matrix representation is $\transferopI : \rr^{N_{\DGSymb}} \rightarrow \rr^{N_{\VelSymb}} $.
\end{itemize} 
To realize the operator $\operator{I}$ we consider the equivalent $L^2$ problem for $u \in \VelSpace$.
\begin{equation}\label{eq:l2I}
\int_{\Omega} (\mathcal{I} \VelVar) \ \Test\, dx = \int_{\Omega} \HdivVar \ \Test\, dx \quad \forall\ \Test \in \DGSpace.
\end{equation}
In terms of coefficient vectors this reads as
\begin{equation}
  \operator{M}^{\DGSymb} (\operator{I} \VelVec) = \operator{M}^{\VelSymb\!,\!\DGSymb} \VelVec, \ \  \forall\ \VelVec \in \rr^{N_{\VelSymb}} \quad \Longrightarrow \quad  \operator{I} = (\operator{M}^{\DGSymb})^{-1} \operator{M}^{\VelSymb\!,\!\DGSymb}
\end{equation}
with the mixed mass matrix 
$$
\operator{M}^{\VelSymb\!,\!\DGSymb}_{i,j} = \int_{\Omega} \varphi_j^{\HdivSymb} \varphi_i^{\DGSymb}\, dx, \ i=1,\hdots,N_{\DGSymb}, j=1,\hdots,N_{\HdivSymb} \quad \text{and} \quad \operator{M}^{\VelSymb\!,\!\DGSymb}_{i,j}=0, j>N_{\HdivSymb}.
$$
Note that $\operator{M}^{\DGSymb}$ is diagonal (for affine linear transformations) such that $(\operator{M}^{\DGSymb})^{-1}$ can be evaluated very efficiently. The overall cost of the transfer operator $\operator{I}$ is essentially that of one sparse matrix multiplication.

Let $\operator{C}^{\DGSymb}: \rr^{N_{\VelSymb}} \times \rr^{N_{\DGSymb}} \rightarrow \rr^{N_{\DGSymb}}$ denote the discrete convection operator corresponding to the trilinearform $\bilinearform{C}$ in \eqref{convstab}.

With the operator $\operator{I}$ we can formulate applications of the convection and the mass operations $\operator{C}^{\DGSymb}(\VelVec), \operator{M}^{\DGSymb}: \rr^{N_{\DGSymb}} \rightarrow \rr^{N_{\DGSymb}}$ with respect to functions in the HDG space $\VelSpace$ and denote the corresponding operators by $\operator{C}^{\VelSymb}(\VelVec), \operator{M}^{\VelSymb}: \rr^{N_{\VelSymb}} \rightarrow \rr^{N_{\VelSymb}}$,
\begin{equation}
  \operator{C}^{\VelSymb}(\VelVec) := \transferopI \operator{C}^{\DGSymb}(\VelVec) \transferop, \qquad \qquad \operator{M}^{\VelSymb} := \transferopI \operator{M}^{\DGSymb} \transferop.
\end{equation}

\begin{remark}[Restriction on time integration scheme] \label{rem:restrtimeint}
Note that the restriction to these two transfer operations implies that we do not allow to apply any part of the Stokes operator to a function in $\DGSpace$ and that no functional on $\VelSpace$ can be used in solution steps involving the convection operator. This is a restriction on the time integration scheme.
\end{remark}

\begin{remark}[Curved elements] \label{rem:embcurv}
If curved elements are considered we no longer have $\HdivSpace \subset \DGSpace$ due to the Piola transform usually applied to construct $H(\divergence)$-conforming finite elements. In this case $\mathcal{I}: U_h \rightarrow V_h$ is not an embedding, but the $L^2$ projection. Nevertheless, $\operator{M}^{\DGSymb}$ is block diagonal with blocks which are not diagonal matrices only on curved elements. Hence, applications of $(\operator{M}^{\DGSymb})^{-1}$ are still cheap.
\end{remark}

\subsection{The semidiscrete formulation}\label{hdgnse:semidiscrete}
With the definitions of the bilinear forms in the previous sections we arrive at the following spatially discrete DAE problem:  
Find $\VelVar(t) \in \VelSpace$ and $p(t) \in \PressureSpace$, such that
\begin{eqnarray} 
%\mbox{ Find } \VelVar(t) \in \VelSpace\mbox{ and } p(t) \in \PressureSpace
%\mbox{, s.t. }\qquad \qquad \qquad \qquad \nonumber \\ 
\label{eq:semidiscrete}
\displaystyle
\left\{
\begin{array}{c@{\hspace*{0cm}}c@{\hspace*{0.05cm}}c@{\hspace*{0.05cm}}cll}
\bilinearform{M}^{\VelSymb}( \frac{\partial}{\partial t}\VelVar,{v}) +
\bilinearform{A}(\VelVar,{v}) + \bilinearform{C}({\HdivVar};\VelVar,{v}) +&
\bilinearform{D}({v},p) 	& = &
\langle {f},v \rangle 	& \forall \ v \in \VelSpace, & t\in [0,T], \\
  & \bilinearform{D}(\VelVar,q)				& = & 0		&
\forall \ q \in \PressureSpace, &t\in [0,T], \\
 &\bilinearform{M} (\VelVar,\VelTest) & = & \bilinearform{M} (\VelVar_0,\VelTest) & \forall \ v \in \VelSpace,		& t=0.
\end{array}
\label{eq:semidiscretense}
\right.
\end{eqnarray}
Here, we implicitly used the embedding $\mathcal{I}$ to define $\bilinearform{C}(\HdivVar;\cdot,\cdot)$ on $\VelSpace$.
Due to \eqref{convstab} we have with 
\begin{equation}
( \frac{\partial}{\partial t}\VelVar,\VelVar)_{L^2} = \frac12 \frac{\mathrm{d}}{\mathrm{d} t} \Vert{\VelVar}\Vert_{L_2}^2 = \Vert{\VelVar}\Vert_{L_2} \frac{\mathrm{d}}{\mathrm{d} t} \Vert{\VelVar}\Vert_{L_2} = \langle f, \VelVar \rangle - \underbrace{\bilinearform{A}(u,u)}_{\geq 0} - \underbrace{\bilinearform{C}({\HdivVar};\VelVar,\VelVar)}_{\geq 0} - \underbrace{\bilinearform{D}(u,p)}_{=0} \leq \Vert f \Vert_{L^2} \Vert \VelVar \Vert_{L^2}
\end{equation}
the stability of the kinetic energy:
\begin{equation}
\frac{\mathrm{d}}{\mathrm{d} t} \Vert{\VelVar}\Vert_{L_2}  \leq \Vert {f}(t) \Vert_{L_2}
\end{equation}
% \cite{dgincnse}\cite{diploma}\cite{ldgincnse}
%
%
%
%
%\section{Tuning stuff}
In the next section we discuss operator splitting time integration methods to solve \eqref{eq:semidiscrete} efficiently.

\section{Operator-splitting time integration} \label{sec:opsplitting}
% \subsection{Introduction to operator splitting approaches}
In this section we are faced with the problem of solving the semi-discrete Navier-Stokes problem with a proper time integration scheme. 
For ease of presentation we neglect external forces ($\operator{f} = 0$) in the following and consider the problem in terms of discrete operators corresponding to the bilinear (trilinear) forms introduced in the previous section: 
Find $\VelVec(t) \in \rr^{N_{\VelSymb}}$ and $\PressureVec(t) \in \rr^{N_{\PressureSymb}}$, such that 
\begin{equation} \label{eq:semidiscop}
\displaystyle
 \left\{
\begin{array}{c @{\hspace{0.1cm}} r @{\hspace{0.1cm}} c @{\hspace{0.1cm}} l @{\hspace{0.1cm}} l}
\operator{M}^{\VelSymb} \frac{\partial \VelVec}{\partial t} + \operator{A} \VelVec 
+ \operator{C}^{\VelSymb}(\VelVec)\ \VelVec  & +\Delta t\, \operator{D}\, \PressureVec  & = & \operator{f} & \text{ in }\ [0,T], \\ 
& \operator{D}^T {\VelVec} & = & 0 & \text{ in }\ [0,T],\\
& \VelVec(t\!=\!0) & = & \VelVec_0. & 
\end{array}
\right.
\end{equation}
In the time integration scheme we want to explicitly exploit the properties of the spatial discretization, i.e. the convection operator $\operator{C}$ should only be involved explicitly in terms of operator evaluation. 
%On the other hand, at the end of every time step we want to ensure the incompressibility constraint. 
Due to the DAE-structure we require time integration schemes which are \emph{stiffly accurate}. Hence, the solution of a Stokes-Brinkman problem should conclude every time step so that the incompressibility constraint is ensured. 
For operator splittings of this kind different approaches exist. We briefly discuss three approaches. 
In section \ref{opsplitting:imex} we discuss \emph{additive} decomposition methods, like the famous class of IMplicit EXplicit (IMEX) schemes. \emph{Product} decomposition methods, sometimes also called \emph{exponential factor} splittings, are discussed in section \ref{opsplitting:oifs} in the framework of operator-integration-factor splittings introduced in \cite{maday1990operator}. 
We discuss the advantages and disadvantages of both approaches and motivate the consideration of a different approach. An operator-splitting modification of the famous fractional step method, cf. \cite{glowinski2003finite}, which eliminates the most important disadvantages of the additive and multiplicative decomposition methods is then introduced and discussed in section \ref{opsplitting:glowinski}. 
We want to stress that the considered operator splitting approaches are of \emph{convection-diffusion} type and should not be confused with projection methods like the Chorin splitting \cite{chorin1967numerical}.

\subsection{Additive decomposition using IMEX schemes}\label{opsplitting:imex}
Additive decomposition methods distinguish spatial operators that are treated implicitly and those that are treated explicitly. The decomposition is \emph{additive} in the sense that every solution (sub-)step involves both operators. 
% the spatial operator in the sense that only one part is treated implicitely and the other part explicitely. A special case of additive decomposition methods are methods which use the same decomposition for all time steps, s.t. one part is only treated explicitely (which is the convection here), i.e. evaluated at old time steps and the other parts also at the new time step. 
This is used by IMplicit EXplicit (IMEX) schemes (see \cite{ascher1995implicit,ascher1997implicit,kanevsky2007application}). The simplest of these schemes is the \emph{semi-implicit Euler} method for which one time step of size $\Delta t$ reads as:
\begin{equation}
  \displaystyle
  \left\{
    \begin{array}{c @{\hspace{0.1cm}} r @{\hspace{0.1cm}} c @{\hspace{0.1cm}} l}
      (\operator{M}^{\VelSymb} + \Delta t \operator{A}) {\VelVec}^{n+1} + & \Delta t \operator{D} \PressureVec^{n+1}  & = & \operator{M}^{\VelSymb} {\VelVec}^{n} - \Delta t \operator{C}^{\VelSymb}(\VelVec^n) \VelVec^n \\ & \operator{D}^T {\VelVec}^{n+1} & = & 0 
    \end{array}
  \right.
  \label{eq:febense}
\end{equation}
Convection only appears on the r.h.s. so that only operator evaluations occur for the convection operator, while the remainder appears fully implicit. This scheme is obviously only first order accurate and only conditionally stable. 
For higher order methods of this decomposition type (e.g. using Multistep or partitioned Runge-Kutta schemes) we refer to \cite{ascher1995implicit,ascher1997implicit,kanevsky2007application}. 
%In the numerical example in section \ref{sec:numex2d} we also used a second order scheme which combines second order implicit and explicit Runge-Kutta methods which share the same time stages, cf. \cite[Table 3.3]{lehrenfeld2010hybrid}.  

The major drawback of this type of operator splitting methods is the fact that the time step size of the explicit and the implicit part of the decomposition have to coincide.
Thereby the stability restriction caused by the explicit treatment of the convection part dictates not only the number of explicit evaluations but also - which is typically more expensive - the number of solution steps for the implicit part. The decomposition methods considered in the sections \ref{opsplitting:oifs} and \ref{opsplitting:glowinski} overcome this issue.

\subsection{Product decomposition with operator-integration-factor splitting }\label{opsplitting:oifs}
 The disadvantage of IMEX methods can be avoided with product decomposition methods, where sequences of 
separated problems are solved successively. The separated problems then only involve the Stokes \emph{or} the convection operator at the same time. The major benefit of this is the fact, that the numerical solution of the sub-problems (e.g. the size of the time steps) can be chosen completely different. 
The derivation of the so called \emph{operator-integration-factor splitting} approach has been derived in \cite[Section 2.1]{maday1990operator}. We use this idea to obtain an operator decoupling for the DAE which allows to \emph{formally} rewrite \eqref{eq:semidiscop} as 
%\vspace*{2cm}
\begin{equation}\label{eq:oifsform}
\left\{
 \begin{array}{rll}
   \frac{\partial}{\partial t} ( Q^{t\rightarrow t^*} \VelVec ) + Q^{t\rightarrow t^*} \operator{M}^{-1} (\operator{A} \VelVec + \operator{D} \PressureVec) &= 0 & \text{ in } [0,T],\\
   \operator{D}^T \VelVec &= 0 & \text{ in } [0,T], 
 \end{array}
\right.
\end{equation}
for some arbitrary $t^* \in \mathbb{R}$ with the integration factor $Q^{t\rightarrow t^*}$, the \emph{propagation operator}, specified later in section \ref{sec:propagate}.

% The step from \eqref{eq:febense} to \eqref{eq:oifsform} is obtained with the following properties of the integration factor:
% \begin{align}
% \frac{\partial}{\partial t} Q^{t \rightarrow t^\ast} = Q^{t \rightarrow t^\ast} \operator{M}^{-1} \operator{C}(\VelVec(t)), \quad Q^{t^\ast \rightarrow t^\ast} = I
% \end{align}

% s.t. $Q^{t_1\rightarrow t_2} \VelVec_1 = \VelVecExpl(t_2)$ with $\VelVecExpl$ the solution of the following \emph{propagation problem}
% \begin{equation}
%  \operator{M} \frac{\partial}{\partial s} \VelVecExpl(s) = -\operator{C}(s) \VelVecExpl(s) \qquad \VelVecExpl(t_1) = \VelVec_1 \label{eq:explproblem}
% \end{equation}
% \vspace*{2cm}
% Note that the identities $Q^{t_b \rightarrow t_c} Q^{t_a \rightarrow t_b} = Q^{t_a \rightarrow t_c}$ and $Q^{t_a \rightarrow t_a} = I$ for $t_a,t_b,t_c \in \mathbb{R}$ hold. 
% We did not specify the space for $\VelVecExpl$ or the operators $\operator{M}$ and $\operator{C}$ in \eqref{eq:oifsform} and \eqref{eq:explproblem}. 
%
Be aware that $\operator{M}^{-1}$ is to be understood only formally for now. 
We are able to apply any suitable (i.e. implicit and stiffly accurate) time integration method for \eqref{eq:oifsform}. 
We do this for a first order method here to explain the procedure and refer to the literature \cite{maday1990operator} for higher order variants.
%A correspondingly defined second order scheme has been used in the numerical examples in section \ref{sec:numex2d}. 

\subsubsection{Implicit Euler.}
We employ the implicit Euler method on \eqref{eq:oifsform} and arrive at
\begin{equation}\label{eq:oifsformiegen}
\left\{
 \begin{array}{rl}
   \frac{1}{\Delta t} ( Q^{t^{n+1}\rightarrow t^*} \VelVec^{n+1} - Q^{t^{n}\rightarrow t^*} \VelVec^{n} ) + Q^{t^{n+1}\rightarrow t^*} \operator{M}^{-1} (\operator{A} \VelVec^{n+1} + \operator{D} \PressureVec^{n+1}) &= 0, \\
   \operator{D}^T \VelVec^{n+1} &= 0.
 \end{array}
\right.
\end{equation}
After setting $t^\ast = t^{n+1}$ and multiplication with $\operator{M}$ this simplifies to
\begin{equation}\label{eq:oifsformie}
\left\{
 \begin{array}{rl}
 \left( \operator{M} + \Delta t \operator{A} \right) \VelVec^{n+1} + \Delta t \operator{D} \PressureVec^{n+1}) &=  \operator{M}^{\DGSymb\!,\VelSymb} Q^{t^{n}\rightarrow t^{n+1}} \transferop \VelVec^{n},  \\
   \operator{D}^T \VelVec^{n+1} &= 0,
 \end{array}
\right.
\end{equation}
which is the solution of a Stokes-Brinkman problem as in \eqref{eq:febense}, but with a different right hand side. 
\subsubsection{The propagation operator.}\label{sec:propagate}
The propagation operator $Q^{t^{1}\rightarrow t^{2}}$ is defined as $Q^{t^{1}\rightarrow t^{2}} \DGVec = \VelVecExpl(t_2)$ where $\VelVecExpl(s)$ solves
\begin{equation} \label{eq:explproblem}
  \frac{\partial \VelVecExpl}{\partial s} + (\operator{M}^{\DGSymb})^{-1} \operator{C}^{\DGSymb}(s) \VelVecExpl(s) = 0, \quad \forall\ s \in (t^1,t^2],\quad \VelVecExpl(t_1) = \DGVec.
\end{equation}
Here $\operator{C}^{\DGSymb}(s) := \operator{C}^{\DGSymb}(\acute{\VelVec}(s))$ with $\acute{\VelVec}(s)$ an \emph{extrapolation} of divergence-free solutions of previous time steps. Note, that the extrapolation of convection velocities renders the problem \eqref{eq:explproblem} \emph{linear} hyperbolic and ensures stability in the sense of \eqref{convstab}.
After replacing $Q$ by a numerical time integrator $Q_{\Delta t}$ for \eqref{eq:explproblem}, the time integration method is completely specified. The order of accuracy of the extrapolation $\acute{\VelVec}$ and the time integrator $Q_{\Delta t}$ should coincide with the order of the time integration method applied on \eqref{eq:oifsform}.
\subsubsection{Properties of product decompositions.}
At this point the advantage of product decomposition methods over IMEX schemes is evident: As the time integration of \eqref{eq:explproblem} is independent of the one for \eqref{eq:oifsform} the (not necessarily) explicit time integrator can deal with stability restrictions (typically using multiple time steps) without influencing the time step for \eqref{eq:oifsform}. This separation of the two problems also introduces the biggest disadvantage of product decomposition methods: an additional consistency error. Even if the DAE \eqref{eq:semidiscop} has a stable stationary solution, product decomposition methods may not reach it. This is not the case for monolithic or additive decomposition approaches. Nevertheless, this \emph{splitting error} is controlled by the time discretization error. 

\begin{remark}[(Marchuk-)Yanenko splitting]
If the implicit Euler discretization as in \eqref{eq:oifsformie} is combined with an Euler method for \eqref{eq:explproblem}, the famous Yanenko splitting method is recovered. 
\end{remark}

\subsection{A modified fractional-step-$\theta$-scheme}\label{opsplitting:glowinski}
In this section we introduce an approach to circumvent severe CFL-restrictions and splitting errors at the same time. We no longer ask for sub-problems that only involve the Stokes or the convection operator (as in section \ref{opsplitting:oifs}), but ask for sub-problems which \emph{only involve one of both implicitly}.
The method is based on \cite{chrispell2007fractional} where the well-known fractional-step-$\theta$-scheme, cf. \cite[Chapter II, section 10]{glowinski2003finite}, is modified. The resulting method is an additive decomposition method without the time step restrictions of IMEX schemes.

We start with \emph{formally} writing down an operator-splitting version of the fractional-step-$\theta$-scheme. The scheme is divided into three steps, the first and the last step treat the Stokes part implicitly and the convection part explicitly as in \eqref{eq:febense} while the second step treats \emph{convection implicitly} and viscosity forces explicitly. 
\begin{subequations} \label{eq:glow} 
\begin{align} \label{eq:glow1} 
\text{Step 1} ( t^{\glowa}\rightarrow t^{\glowb} ):&  \left\{
    \begin{array}{c @{\hspace{0.1cm}} r @{\hspace{0.1cm}} c @{\hspace{0.1cm}} l}
      (\operator{M} + \theta \Delta t \operator{A}) {\VelVec}^{\glowb} + & \theta \Delta t \, \operator{D} \, \PressureVec^{\glowb}  & = & \operator{M} {\VelVec}^{\glowa} - \theta \Delta t \operator{C} \VelVec^{\glowa} \\ & \operator{D}^T {\VelVec}^{\glowb} & = & 0 
    \end{array}
  \right. \\ \label{eq:glow2} 
\text{Step 2 } ( t^{\glowb}\rightarrow t^{\glowc} ):&  \left\{ 
\hspace*{0.23cm}
    (\operator{M} + \theta^\ast \Delta t \operator{C}) \VelVec^{\glowc} 
\hspace*{1.74cm} = 
\operator{M} {\VelVec}^{\glowb} - \theta^\ast \Delta t (\operator{A} \VelVec^{\glowb} + \operator{D} \PressureVec^{\glowb})
  \right. \\ \label{eq:glow3} 
\text{Step 3 } ( t^{\glowc}\rightarrow t^{\glowd} ):&  \left\{
    \begin{array}{c @{\hspace{0.1cm}} r @{\hspace{0.1cm}} c @{\hspace{0.1cm}} l}
      (\operator{M} + \theta \Delta t \operator{A}) {\VelVec}^{\glowd} + & \theta \Delta t \, \operator{D} \,  \PressureVec^{\glowd}  & = & \operator{M} {\VelVec}^{\glowc} - \theta \Delta t \operator{C} \VelVec^{\glowc} \\ & \operator{D}^T {\VelVec}^{\glowd} & = & 0 
    \end{array}
  \right.
\end{align}
\end{subequations}
The time stages are labeled by the superscripts $\glowa$, $\glowb$, $\glowc$, $\glowd$,  respectively, where $\glowa$ denotes initial data and $\glowd$ the final time stage. The time steps size for the first and the last step is $\theta \Delta t$ and $\theta^\ast \Delta t$ in the middle step with $\theta^\ast = 1-2 \theta$. Here $\theta = 1 - 1/\sqrt{2}$. 
We note that specifications for $\operator{M}$ and $\operator{C}$ with respect to the considered spaces and the convection velocity for $\operator{C}$ are still missing at this points. 
The scheme is second order accurate. In contrast to the unsplit fractional-step-$\theta$-scheme the stability analysis of this time integration scheme is an open problem. In our experience, however, the stability restrictions are much less restrictive than those of a comparable IMEX schemes.

% \begin{remark}[Stability]
% The stability of the scheme in \eqref{eq:glow} is not unconditionally stable. \todo[inline]{Which kind of statements to make here??!}
% \end{remark}
\subsubsection{Sub-steps in different spaces.}
Initially, we stated that we want to avoid solving linear systems involving convection, for efficiency reasons. This seems to be contradictory to what is formulated in \eqref{eq:glow2}. Nevertheless, as only convection is involved implicitly an efficient numerical solution is still possible. We apply a simple iterative scheme which only involves explicit operator evaluations of the convection to do so. We explain this in section \ref{sec:pseudotime}. To do this efficiently, we want Step 2 to be formulated in the space $\DGSpace$ while Step 1 and 3 are to be formulated in the space $(\VelSpace,\PressureSpace)$. This poses problems the solution of which we discuss in this section.
 
In Step 1 and Step 3 the adjustments are obvious: 
Step 1 only depends on initial data in $\VelSpace$. The initial data for Step 3 is in $\DGSpace$ but appears only in terms of functionals ($\operator{M} \VelVecExpl$ and $\operator{C} \VelVecExpl$) such that the transfer operations are clear.  
Step 2 is more involved, cf. remark \ref{rem:restrtimeint}. Functionals in $\VelSpace'$ (such as $\operator{A} \VelVec$) are in general not functionals in $\DGSpace'$. We use the first equation in \eqref{eq:glow1} to formally define a different representation of the functionals required in Step 2:
\begin{equation}
\theta \Delta t \ \operator{g}^{\glowb} = - \theta \Delta t (\operator{A} {\VelVec}^{\glowb} + \operator{D} \PressureVec^{\glowb})  = \operator{M} ({\VelVec}^{\glowb} - {\VelVec}^{\glowa})  + \theta \Delta t \operator{C} \VelVec^{\glowa}
\end{equation}
Now we replace $\operator{M}$ and $\operator{C}$ with operations suitable for a setting in $\DGSpace$:
\begin{equation}\label{eq:g}
\rr^{N_{\DGSymb}} \ni \operator{g}^{\glowb} := \frac{1}{\theta \Delta t} \ \operator{M}^{\VelSymb\!,\!\DGSymb} ({\VelVec}^{\glowb} - {\VelVec}^{\glowa})  + \theta \Delta t \operator{C}^{\DGSymb}(\VelVec^{\glowa}) \transferop \VelVec^{\glowa}
\end{equation}
We arrive at the modified fractional-step-$\theta$-scheme:
\begin{subequations} \label{eq:modglow} 
\begin{align} \label{eq:modglow1} 
\text{Step 1}:&  \left\{
    \begin{array}{c @{\hspace{0.1cm}} r @{\hspace{0.1cm}} c @{\hspace{0.1cm}} l}
 % \hspace*{0.278cm}
      (\operator{M}^{\VelSymb}\!\! + \theta \Delta t \operator{A}) {\VelVec}^{\glowb} + & \theta \Delta t \operator{D}\ \PressureVec^{\glowb}  & = & \operator{M}^{\VelSymb} {\VelVec}^{\glowa} - \theta \Delta t \ \transferopI \operator{C}^{\DGSymb}\! (\VelVec^{\glowa}) \ \transferop \ \VelVec^{\glowa} \\ & \operator{D}^T {\VelVec}^{\glowb} & = & 0 
    \end{array}
  \right. \\ \label{eq:modglow2} 
\text{Step 2}:&  \left\{ 
 \hspace*{0.25cm}
    (\operator{M}^{\DGSymb}\!\! + \theta^\ast \Delta t \operator{C}^{\DGSymb}\!(t^{\glowc}) \VelVecExpl^{\glowc} 
\hspace*{0.985cm} 
= 
\operator{M}^{\VelSymb,\DGSymb}\!\! {\VelVec}^{\glowb} - \theta^\ast \Delta t \operator{g}^{\glowb}
  \right. \\ \label{eq:modglow3} 
\text{Step 3}:&  \left\{
    \begin{array}{c @{\hspace{0.1cm}} r @{\hspace{0.1cm}} c @{\hspace{0.1cm}} l}
 % \hspace*{0.278cm}
      (\operator{M}^{\VelSymb}\!\! + \theta \Delta t \operator{A}) {\VelVec}^{\glowd} + & \theta \Delta t \operator{D} \ \PressureVec^{\glowd}  & = & \operator{M}^{\DGSymb,\VelSymb} {\VelVecExpl}^{\glowc} \!\!- \theta \Delta t \ \transferopI \operator{C}^{\DGSymb}\!(t^{\glowc})  \VelVecExpl^{\glowc} \\ & \operator{D}^T {\VelVec}^{\glowd} & = & 0 
    \end{array}
  \right.
\end{align}
\end{subequations}
where we replaced the generic operators $\operator{M}$, $\operator{C}$ with suitable ones. We recall the definition of the extrapolated convection operator $\operator{C}^{\DGSymb}\!(t^{\glowc}) := \operator{C}^{\DGSymb}\!(\acute{\VelVec}(t^{\glowc}))$ the convection velocity of which is (linearly) extrapolated from exactly divergence-free velocities. 
% We first discuss the easier steps, Step 1 and 3. 
% With $\operator{M} = \operator{M}^{\VelSymb}$ and $\operator{C} = \operator{C}^{\VelSymb}(\VelVec^n)$ in \eqref{eq:glow1} Step 1 can be defined adequately. 
% In Step 3 we set $\operator{M} = \operator{M}^{\VelSymb}$ on the l.h.s. and  $\operator{M} \VelVec^{\glowc} = \operator{M}^{\DGSymb,\VelSymb} \VelVecExpl^{\glowc}$ on the r.h.s.. Further we set $\operator{C} \VelVec^{\glowc} = \transferopI \operator{C}^{\DGSymb}(\acute{\VelVec}(t^{\glowc})) \VelVecExpl^{\glowc}$ with the extrapolated (divergence-free) velocity $\acute{\VelVec}$ as in section \ref{sec:propagate}.

% First, we need to apply the transfer operations $\transferop$ to switch between solutions in $\DGSpace$ and $(\VelSpace)$ . This allows us to carry out Step 1 and Step 3 if $\VelVec^{n} \in \rr^{N_{\VelSymb}}$ and $\VelVec^{\glowc} \in \rr^{N_{\DGSymb}}$.
% Secondly, 

\subsubsection{Iterative solution of the implicit convection problem.}\label{sec:pseudotime}
In \eqref{eq:modglow2} we need to solve a problem of the form 
\begin{equation} \label{eq:implconv}
  \VelVecExpl + \tau^\ast (\operator{M}^{\DGSymb})^{-1} \operator{C}^{\DGSymb} \VelVecExpl = \operator{g}^\ast
\end{equation}
for given $\tau^\ast$, $\operator{g}^\ast$ and constant convection $\operator{C}^{\DGSymb}$. We do this by means of a pseudo time-stepping method, i.e. we formulate \eqref{eq:implconv} as the stationary solution to
\begin{equation}
  \frac{\partial}{\partial s} \VelVecExpl(s) + \VelVecExpl(s) + \tau^\ast \operator{C}^{\DGSymb} \VelVecExpl(s) = \operator{g}^\ast, \quad s \in [0,\infty).
\end{equation}
Note that a stationary solution exists as $\mathbf{Id} + \tau^\ast \operator{C}^{\DGSymb}$ only has eigenvalues $\lambda$, with $Re(\lambda) > 1$. This stationary solution is approximated with a few explicit Euler time step with an artificial time step size $\Delta s$. 
\begin{equation} \label{eq:implconvit}
 \VelVecExpl_{i+1} =  \VelVecExpl_{i} + \Delta s (\operator{g}^\ast - \VelVecExpl_i + \tau^\ast 
(\operator{M}^{\DGSymb})^{-1} \operator{C}^{\DGSymb} \VelVecExpl_i).
\end{equation}
This procedure can be interpreted as a Richardson iteration applied to \eqref{eq:implconv}. 
With a time step size which is tailored for stability (as in the numerical solution to \eqref{eq:explproblem}) we iterate \eqref{eq:implconvit} until the initial residual is reduced by a prescribed factor. 
To solve the convection step, Step 2, we only require operator evaluations for the convection. The time step size $\Delta s$ in this iteration is decoupled from the time step size $\Delta t$ used of the overall scheme.
%As an initial guess we take $\VelVecExpl_0 = \operator{g}^\ast$ which coincides with an explicit Euler step of size $\tau^\ast$ to 

%\subsubsection{Properties of the modified fractional-step-$\theta$-scheme}
% \todo[inline]{Clarify pros and cons of operator slitting methods}
\section{Numerical examples} \label{sec:numex}
In this section we consider different test problems which essentially purpose three different goals: 
\begin{enumerate}
\item The validation of convergence properties of the HDG method with and without the projected jumps.
\item The investigation and quantification of the dependency of the sparsity pattern of arising linear systems on the choice of the discrete velocity spaces.
\item The evaluation of the performance of the space and time discretization on benchmark problem.
\end{enumerate}
The examples in section \ref{sec:potsol} and section \ref{sec:kovasznay} aim at goal 1 and are two-dimensional test cases with known exact solutions for a stationary Stokes and Navier-Stokes problem, respectively. These cases allow for a thorough investigation of the convergence history of the proposed spatial discretizations in the usual norms. 
The test case in section \ref{sec:numexlinsys} aims at goal 2 and
is concerned with the impact of the choice of discretization spaces on the complexity of linear systems. For this purpose, a three-dimensional vector-valued Poisson problem is considered using different $H(\divergence)$-conforming DG and HDG spaces. The impact of hybridization and the improvement due to the projected jumps modification is compared to other DG methods. 
Finally, we approach goal 3 by considering two challenging transient benchmark problems in sections \ref{sec:numex2d} and \ref{sec:numexbenchm}. The benchmark problems have been defined within the DFG Priority Research Programm 'Flow Simulation on High Performance Computers' and are formulated in \cite{schafer1996benchmark}.
In section \ref{sec:numex2d} a two-dimensional benchmark problem is considered to demonstrate the accuracy of our method for a demanding test case and to compare the discussed time integration methods. 
Finally, in section \ref{sec:numexbenchm} we discuss a three-dimensional, and hence computationally demanding, benchmark problem from \cite{schafer1996benchmark}. We compare accuracy and run-time performance to the data of the studies in \cite{schafer1996benchmark,bayraktar2012benchmark,john2006efficiency}. 

The methods discussed in this paper have been implemented in the add-on package \texttt{ngsflow} \cite{ngsflow} for the high order finite element library \texttt{NGSolve} \cite{schoeberl2014cpp11}. The computations in this section have also been carried out with this software. Throughout this section we only consider direct solvers and comment on linear solvers below, in remark \ref{rem:precond}.

For the computations of the benchmark problem in sections \ref{sec:numex2d} and \ref{sec:numexbenchm} we used the reduced $H(div)$-conforming space $\HdivSpace^{\text{red}}$, cf. remark \ref{rem:red} and the \emph{projected jumps} modification, cf. section \ref{sec:modifications}.

\subsection{Stokes flow around obstacle} \label{sec:potsol}
We consider the Stokes problem in the domain $\Omega=[-2,2]^2 \setminus \Omega^{-}$ with the circular obstacle $\Omega^{-} := \{\Vert x \Vert \leq 1 \}$ and viscosity $\nu=1$. On the whole boundary except for the outflow boundary $\{x=2\}$ we prescribe Dirichlet boundary conditions, on $\{ x=2\}$ we prescribe Neumann-kind boundary conditions such that the solution to the problem is given by $u = (\partial_y \Psi, - \partial_x \Psi)$ with the potential $\Psi = y (1-\frac{1}{x^2+y^2})$. Note that $\Psi|_{\partial \Omega^{-}}$ is constant so that $u \cdot n = \nabla \Psi \times n = 0$. On $\partial \Omega^{-}$ we have $u \cdot n = 0$, but $u \neq 0$, i.e. we have a slip on the obstacle. Further we have that $p = 0$ in $\Omega$.

\begin{figure}[h!]
    \hspace*{-0.5cm}
    \includegraphics[width=0.35\textwidth]{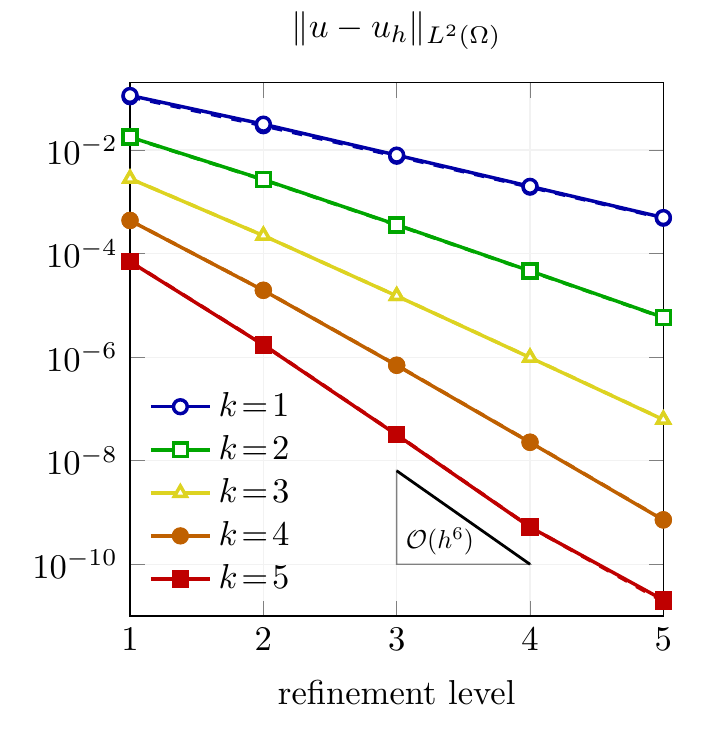}
    \hspace*{-0.5cm}
    \includegraphics[width=0.35\textwidth]{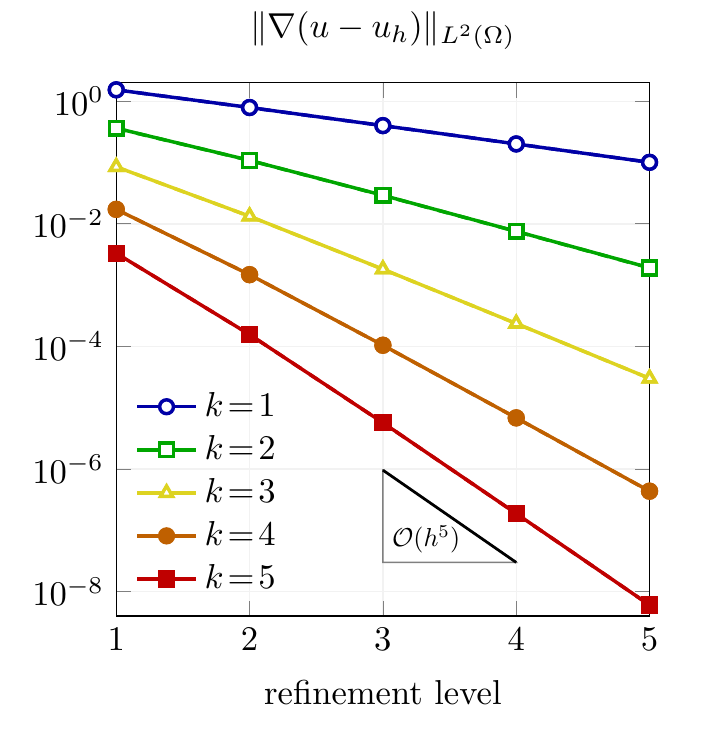}
    \hspace*{-0.5cm}
    \includegraphics[width=0.35\textwidth]{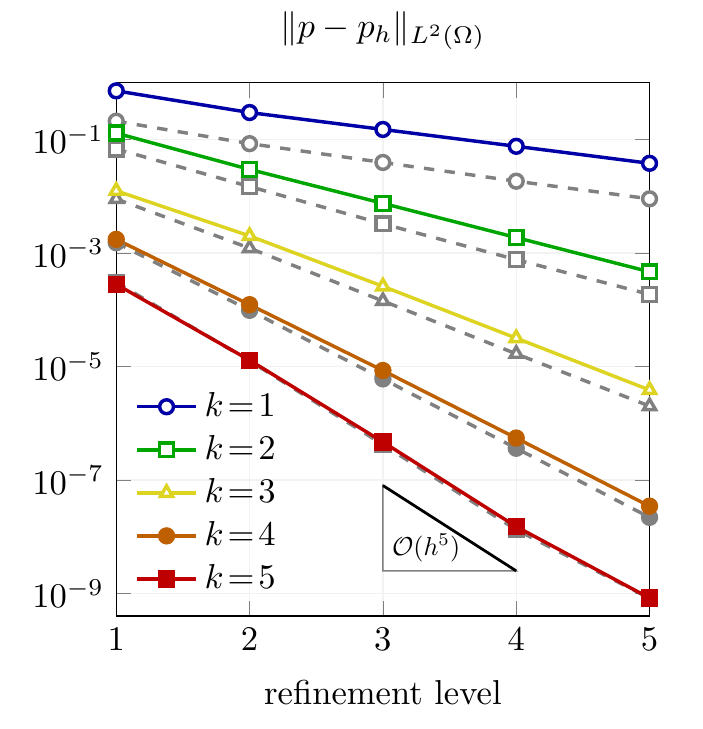}
    \hspace*{-0.5cm}
    \vspace*{-0.2cm}
\caption{Convergence of divergence-conforming HDG method with (solid) and without (dotted) projected jumps in different norms for the example in section \ref{sec:potsol}. The dotted lines in gray correspond to the HDG method without projected jumps, the solid lines in color correspond to the HDG method with projected jumps. Except for the error in the pressure the results are hardly distinguishable.}
\label{fig:potsol}
\end{figure}

On an initially coarse unstructured grid with only 72 triangles we obtain the results displayed in Figure \ref{fig:potsol} after succesive uniform mesh refinements. We compare the exact solution $u$, $p$ and the computed solution $u_h$, $p_h$ and investigate the convergence of the error in the norms $\Vert u - u_h \Vert_{L^2(\Omega)}$, $\Vert \nabla (u - u_h) \Vert_{L^2(\Omega)}$ and $\Vert p - p_h \Vert_{L^2(\Omega)}$.
We consider the discretization with and without the projected jumps modification. In all norms we observe optimal order convergence, i.e.  $\Vert u - u_h \Vert_{L^2(\Omega)} = \mathcal{O}(h^{k+1})$, $\Vert \nabla (u - u_h) \Vert_{L^2(\Omega)} = \mathcal{O}(h^{k}) $ and $\Vert p - p_h \Vert_{L^2(\Omega)} = \mathcal{O}(h^{k})$ where $k$ is the order of the velocity field inside the elements. 
Note that after static condensation the remaining degrees of freedoms are the unknowns corresponding to order $k$ polynomials for the tangential and order $k$ polynomials for the normal component on the facets and one constant per element for the pressure. For the projected jumps formulation the tangential unknowns are reduced by one order. 
The difference between the error of both formulations --- with and without projected jumps --- is only marginal in the velocity. In the pressure field we observe a difference between the methods, however both converge with optimal order and the difference decreases for increasing polynomial degree $k$. 
For the impact of the projected jumps modification concerning the computational effort, we refer to section \ref{sec:numexlinsys} where this aspect is discussed in more detail.

\subsection{Kovasznay flow} \label{sec:kovasznay}
As a test case for the stationary Navier-Stokes equations, we consider the famous example by \cite{kovasznay1948laminar}.
On the boundary of the domain $\Omega = [-\frac12,\frac32] \times [0,2]$ we again prescribe inhomogeneous Dirichlet data, so that the exact solution to the Navier-Stokes equations (with $\nu = 1$) is
\begin{equation*}
u(x,y) = \left( \begin{array}{c} 1 - e^{\lambda x} \cos (2\pi y) \\ \frac{\lambda}{2 \pi} e^{\lambda x} \sin( 2 \pi y) \end{array} \right), \quad p(x,y) = - \frac12 e^{2 \lambda x} + \bar p \text{ with } \bar p \in \mathbb{R}
\end{equation*}
Here, $\lambda = \frac{-8\pi^2}{\nu^{-1} + \sqrt{\nu^{-2} + 64 \pi^2}}$ and $\bar p$ so that $\int_{\Omega} p \, d x = 0$. 

\begin{figure}
    \hspace*{-0.5cm}
    \includegraphics[width=0.35\textwidth]{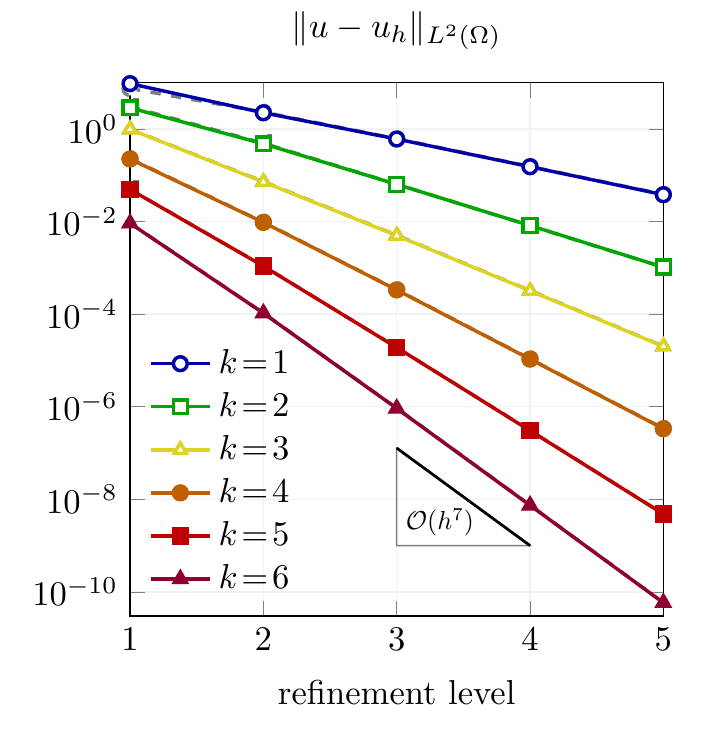}
    \hspace*{-0.5cm}
    \includegraphics[width=0.35\textwidth]{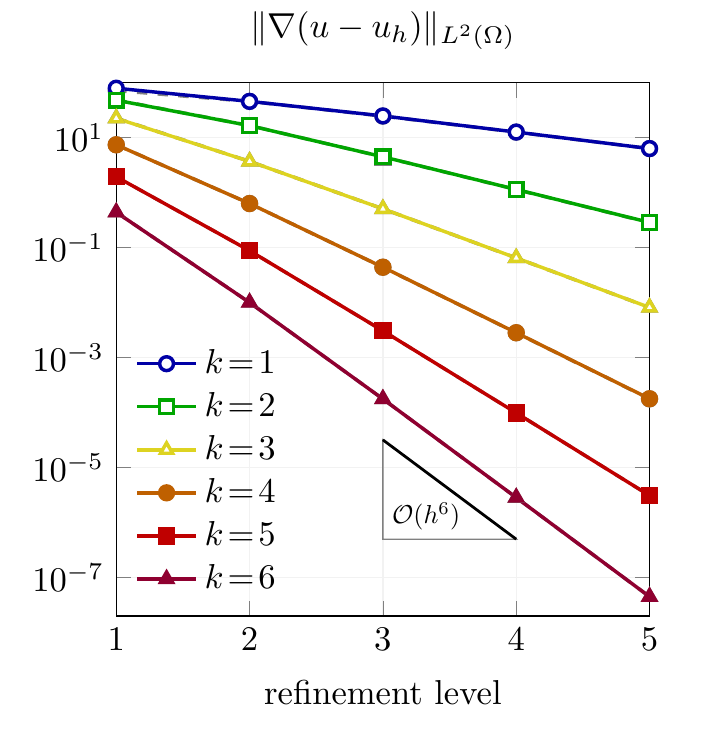}
    \hspace*{-0.5cm}
    \includegraphics[width=0.35\textwidth]{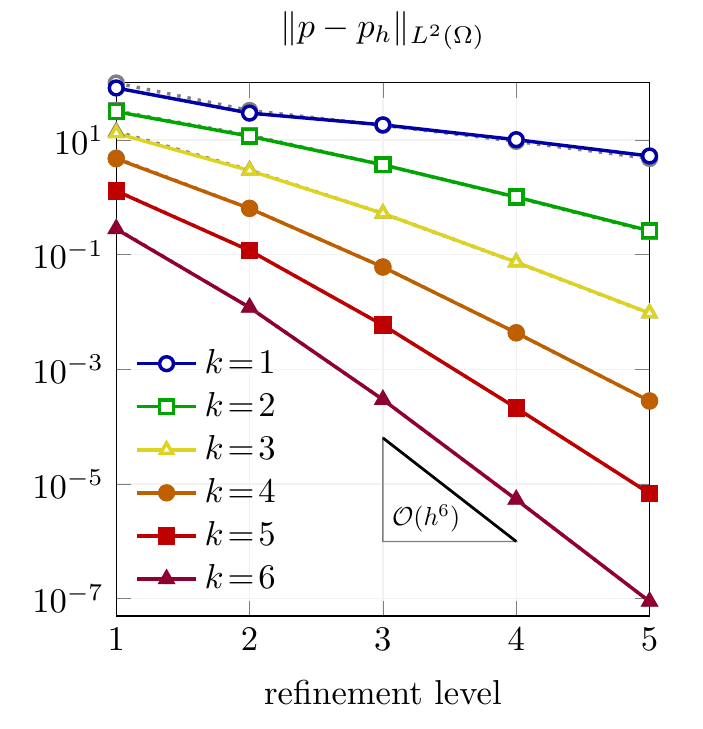}
    \hspace*{-0.5cm}
\caption{Convergence of divergence-conforming HDG method with (solid) and without (dotted) projected jumps in different norms for the example in section \ref{sec:kovasznay}. The results are almost identical.}
\label{fig:kovasznay}
\end{figure}

To approximate the stationary solution we initially solve the Stokes problem corresponding to the boundary data and use the 
simple IMEX scheme of first order, cf. \eqref{eq:febense}, to progress to $t=1000$ with 50 time steps of size $\Delta t = 20$. This choice of time discretization parameters gives stable solutions for all considered spatial discretizations but at the same time provides time discretization errors which are neglegible compared to the spatial errors. 
On an initally coarse unstructured grid with only 18 triangles we obtain the results displayed in Figure \ref{fig:kovasznay} after succesive uniform mesh refinements. We observe that the errors converge optimal in the considered norms for the pressure and the velocity. Further, we observe that the difference between the results obtained with and without the projected jumps formulation, e.g. with and without a reduction of the degrees of freedoms at the facets, is again only marginal. In fact one only observes a (very small) difference in the case $k=1$.

\subsection{Linear systems - A comparison between DG and HDG methods} \label{sec:numexlinsys}
We consider the comparably simple vector-valued Poisson problem
\begin{equation} \label{eq:vecpoisson}
 -\Delta u = f  \text{ in } \Omega, \quad u = 0 \text{ on } \partial \Omega,
\end{equation}
discretized with four different methods.
The three-dimensional domain $\Omega$ is the same as the one in section \ref{sec:numexbenchm} with an unstructured tetrahedral mesh consisting of 3487 elements. The problem \eqref{eq:vecpoisson} leads to discretizations with symmetric positive definite matrices $\operator{A}$. Only the lower triangular part of the sparse matrix has to be stored and a sparse Cholesky factorization algorithm can be used. In this comparison, we used the sparse direct solver \texttt{PARDISO}, cf. \cite{pardiso,schenk04:_solvin_unsym_spars_system_linear_equat_pardis}. 
We are only concerned with the sparsity pattern of different methods before and after static condensation and do not compare accuracy or conditioning of the discretizations. Moreover, for the sake of simplicity we do \emph{not apply} the reduction of the space $\HdivSpace$ as discussed in remark \ref{rem:red}. 
The following quantities of interest for the linear systems arising from discretizations of \eqref{eq:vecpoisson} are displayed in Table \ref{table1}:
\begin{tabular}{p{1.5cm}@{\hspace{0.1cm}:\hspace{0.1cm}}l}
 \#dof[K]& number of unknowns (in thousands)\\
 \#cdof[K]& number of unknowns after static condensation (in thousands)\\
 \#nzeA[K]& number of non-zero entries in the system matrix $\operator{A}$ (lower triangle) (in thousands)\\
 \#nzeL[K]& number of non-zero entries in the Cholesky factor $\operator{L}$ (in thousands)
\end{tabular} 
\vspace*{0.2cm}\\
Four different methods are considered on $\HdivSpace$ or $\VelSpace$ with varying polynomial degree between $k=1$ and $k=6$:
\begin{enumerate}
\item HDG: The HDG method proposed in section \ref{sec:DGHDG} without \emph{projected jumps}.
%: $\bilinearform{A}(\VelVar,\VelTest) = \langle f , \VelTest \rangle, \forall \VelTest \in \VelSpace$. 
\item PHDG: The HDG method with \emph{projected jumps}, cf. section \ref{sec:modifications}.
%: $\bilinearform{A}^r(\VelVar,\VelTest) = \langle f , \VelTest \rangle, \ \forall \ \VelTest \in \VelSpace$.
\item Std.DG.: A standard DG method using the space $\HdivSpace$ where the basis is constructed such that all degrees of freedom from one element couple with \emph{all} degrees of freedom from adjacent elements. 
\item N.DG: A nodal DG method using the space $\HdivSpace$ where basis functions are assumed to be constructed such that degrees of freedom associated to one element couple only with degrees of freedom from adjacent elements which have support on the shared facet. At the same time we assume that basis functions are constructed such that the number of basis functions with support on a facets is minimized, cf. \cite{hesthaven2007nodal}. In terms of the sparsity pattern this nodal DG method represents the best case for a DG method without hybridization.
\end{enumerate}
\begin{table}[h!]
\scriptsize
  \begin{center}
    \begin{tabular}{l@{\ \ \ }r@{\ \ \ }r@{\ \ \ }r@{\ \ \ }r|r@{\ \ \ }r@{\ \ \ }r@{\ \ \ }r}
      \toprule
      & \small Std.DG & \small N.DG &  \small HDG & \small PHDG 
      & \small Std.DG & \small N.DG &  \small HDG & \small PHDG \\
      \midrule
      & & \multicolumn{2}{c}{$k=1$} & &
      & \multicolumn{2}{c}{$k=2$} & \\
      \midrule
      \#dof[K]   &23&23&69&38 & 67 & 67 & 158 & 112 \\
      \#cdof[K]   &23&23&69&38 & 67 & 67 & 137 &  91 \\
      \#nzeA[K]   &732&676&2 037&637 & 5\ 073& 4\ 177& 8\ 103& 3\ 628\\
      \#nzeL[K]   &10\ 113&10\ 208&17 768&5\ 569 & 64\ 463&69\ 831&70\ 138&31\ 415\\
      \midrule
      & & \multicolumn{2}{c}{$k=3$} & &
      & \multicolumn{2}{c}{$k=4$} & \\
      \midrule
      \#dof[K]   & 146 & 146 & 298 & 237 & 271 & 271 & 500 & 423 \\
      \#cdof[K]   & 146 & 146 & 229 & 168 & 271 & 261 & 343 & 267 \\
      \#nzeA[K]   & 21\ 686 & 16\ 051 & 22\ 443 & 12\ 124 &  69\ 525 &  46\ 912 &  50\ 412 & 30\ 588 \\
      \#nzeL[K]   & 261\ 977 & 260\ 416 & 194\ 524 &104\ 496 & 814\ 168 & 731\ 253 & 435\ 129 & 264\ 183 \\
      \midrule
      & & \multicolumn{2}{c}{$k=5$} & &
      & \multicolumn{2}{c}{$k=6$} & \\
      \midrule
      \#dof[K]   & 453 & 453 & 773 & 681 & 702 & 702 & 1\ 128 & 1\ 022\\
      \#cdof[K]   & 453 & 411 & 480 & 389 & 702 & 597 & 640 & 533\\
      \#nzeA[K]   & 184\ 035 & 101\ 473 & 98\ 696 & 64\ 816 & 424\ 764 & 200\ 813 & 175\ 321 & 121\ 944\\
      \#nzeL[K]   & 2\ 099\ 690 & 1\ 741\ 072 & 847\ 913 & 557\ 798 & 4\ 752\ 072 & 3\ 519\ 061 & 1\ 502\ 558 & 1\ 045\ 444\\
      \bottomrule
    \end{tabular}
  \end{center}
\caption{Comparison of different DG methods for the vector-valued reaction Poisson problem \eqref{eq:vecpoisson}.}
\label{table1}
\end{table}
In Table \ref{table1} we observe that for small polynomial degree $k$ the amount of additional unknowns required for the HDG formulation is quite large as are the nonzero entries in the system matrix. Nevertheless except for $k=1$, the number of nonzero entries in the Cholesky factor are comparable ($k=2$) to the Standard DG methods or less ($k>2$). For high order, i.e. $k\geq4$ the HDG method performs significantly better as it has less nonzero entries in $\operator{A}$ and $\operator{L}$. 
The HDG space with the \emph{projected jumps} modification improves the situation dramatically. It essentially compensates the overhead of the HDG method for small $k$. But even for $k=6$ the effect is still significant. We note, that the difference between the first three methods increases for increasing polynomial degree $k$ whereas the difference between the last two method, PHDG and HDG, decreases. In all cases the HDG method with \emph{projected jumps} outperforms all alternatives. 

\subsection{A two-dimensional benchmark problem} \label{sec:numex2d}
In this section we consider the benchmark problem denoted as ``2D-2Z'' in \cite{schafer1996benchmark} where a laminar flow around a circle-shaped obstacle is considered. The Reynolds number is moderately high ($Re = 100$) and results in a periodic vortex street behind the obstacle. 
We briefly introduce the problem (for more details we refer to \cite{schafer1996benchmark}) and the numerical setup to investigate spatial and temporal discretization errors. Finally, we discuss the obtained results.

% \subsubsection{Geometry and boundary condition}
% \paragraph{2D Geometry} \ \\
% \begin{figure}[h]
% \centering
%   \begin{overpic}[width=0.95\textwidth]{../graphics/domain2d.png}
%      \put(2,10){\LARGE $\Gamma_{in}$}
%      \put(335,10){\LARGE $\Gamma_{out}$}
%  \end{overpic}
% %\includegraphics[width=0.95\textwidth]{../graphics/domain2d.png}
% \caption{2D geometry for test cases of section \ref{numexamples:navstokes}}
% \label{fig:domain2D}
% \end{figure}
\subsubsection{Geometrical setup and boundary conditions.}
The domain is a rectangular channel without an almost vertically centered circular obstacle, cf. Figure \ref{fig:sketch2d},
\begin{equation}
  \Omega := [0,2.2] \! \times \! [0,0.41] \setminus \{ \Vert x - (0.2,0.2) \Vert_2 \leq 0.05\}.
\end{equation}
The boundary is decomposed into $\Gamma_{in} := \{x = 0\}$, the inflow boundary, $\Gamma_{out} := \{x = 2.2\}$, the outflow boundary and $\Gamma_W := \partial \Omega \setminus ( \Gamma_{in} \cup \Gamma_{out} )$, the wall boundary. On $\Gamma_{out}$ we prescribe natural boundary conditions $ ( - \nu \nabla u + p I ) \cdot \normal = 0 $, on $\Gamma_W$ homogeneous Dirichlet boundary conditions for the velocity and on $\Gamma_{in}$ the inflow Dirichlet boundary conditions
\begin{equation*}
 u(0,y,t) = u_D = (3/2 \cdot \bar{u}) \cdot 4 \cdot {y (d_y - y)}/{d_y^2} \cdot (1,0,0).
\end{equation*}
Here, $\bar{u} = 1$ and the viscosity is fixed to $\nu = 10^{-3}$ which results in a Reynolds number $Re=100$. 

\subsubsection{Drag and Lift.} 
The quantities of interest in this example are the (maximal and minimal) drag and lift forces $c_D$, $c_L$ that act on the disc. These are defined as
\begin{equation*}
 c_D := \frac{1}{\bar{u}^2 r} \int_{\Gamma_\circ} \left( \nu \frac{\partial u}{\partial \normal} - p \normal \right) \cdot e_x \ ds,
 \quad \quad  
 c_L := \frac{1}{\bar{u}^2 r} \int_{\Gamma_\circ} \left( \nu \frac{\partial u}{\partial \normal} - p \normal \right) \cdot e_y \ ds.
\end{equation*}
Here $e_x, e_y$ denote the unit vectors in $x$ and $y$ direction, $r = 0.05$ is the radius of the obstacle, $\bar{u}$ is the average inflow velocity ($\bar{u}=1$) and $\Gamma_\circ$ denotes the surface of the obstacle. 
%If the flow is periodic in time you can further define the \emph{Strouhal} number $St := \frac{2 r}{\bar{u} T} $ with $T$ the length of one period.
\subsubsection{Numerical setup.} 
\begin{figure}[h!]
  \begin{center}
    \includegraphics[height=0.145\textwidth]{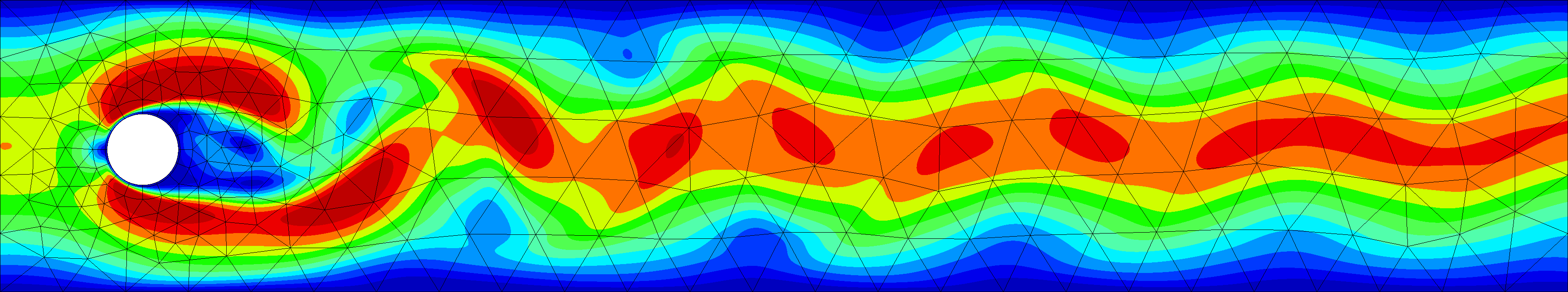}
    \hfill
    \begin{overpic}[width=0.145\textwidth,angle=90]{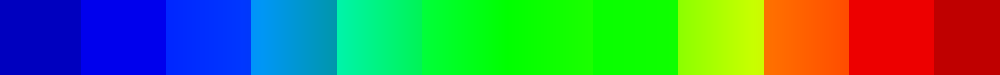}
      \put(7,0){\small $0$}
      \put(7,53){\small $2$}
    \end{overpic}
    \hfill
    \hfill
    \hfill
    \includegraphics[height=0.145\textwidth]{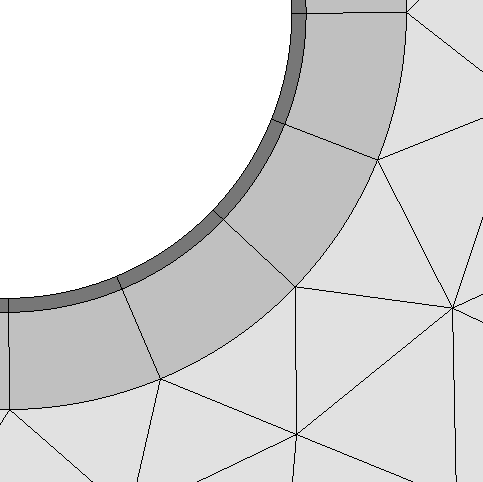}
  \end{center}
  \vspace*{-0.5cm}
  \caption{Sketch of the mesh and the solution (color coding corresponding to velocity magnitude $\Vert \VelVar \Vert_2$) to the problem considered in section \ref{sec:numex2d} at a fixed time $t$ (left) and zoom-in on the boundary layer mesh (right).}
  \label{fig:sketch2d}
\end{figure}

We use an unstructured triangular grid with an additional layer of quadrilaterals around the disk which is anisotropically refined towards the disk once. In Figure \ref{fig:sketch2d} the geometry, the mesh and a typical solution is depicted. 

In order to be able to neglect time discretization errors, when investigating the spatial accuracy,  we consider the use of a well-known stiffly accurate second order Runge-Kutta-IMEX scheme, taken from \cite[section 2.6]{ascher1997implicit}, with an extremely small time step size $3.125 \cdot 10^{-5}$ which means that roughly $10000$ time steps are used to resolve one full period. 
The duration of one full cycle is roughly $1/3 s$.
We also fix the mesh and consider a pure $p$-refinement, i.e. variations of the polynomial degree $k$.

For the investigations of the temporal discretization error, we consider a fixed polynomial degree $k=6$, with the same mesh. For the considered example the stability restriction of the second order IMEX scheme is severe. A time step of below $10^{-3}$ has to be considered. The intention of the discussion of operator-splitting methods in section \ref{sec:opsplitting} has been to present strategies to circumvent or relax these severe conditions. We compare the performance of a product decomposition method, a second order operator-integration-factor splitting version of the BDF2 method, and the modified fractional-step-$\theta$-scheme discussed in section \ref{opsplitting:glowinski}. Note that these methods allows to consider much larger time steps than the IMEX scheme. As references we give the values from the literature \cite{featflow} and from the IMEX scheme with $\Delta t=10^{-3}$ (stability limit) and the reference solution with $\Delta t=3.125 \cdot 10^{-5}$.

\subsubsection{Numerical results: spatial discretization.} 
In Table \ref{tab:2d} the quantities of interest are shown for varying polynomial degree $k$. As a reference we also show the result obtained by \texttt{FEATFLOW} \cite{featflowman} with a discretization using quadrilateral meshes and continuous second order finite elements for the velocity with a discontinuous piecewise linear pressure ($Q_2/P_1^{\text{disc}}$). These results have been made accessible on \cite{featflow}.
%If a digit in the table is underlined, we consider the quantity converged until that digit. 
\begin{table}[h!]
  \begin{center}
\small
% BEGIN RECEIVE ORGTBL tab:2d
\begin{tabular}{crrrrr}
\toprule
 & \#dof & $\max c_D$ & $\min c_D$ & $\max c_L$ & $\min c_L$ \\
% \midrule
% 1 &  1\ 380 & 2.73697 & 2.73697 & -0.15591 & -0.15591 \\
% 2 &  2\ 582 & 3.20792 & 3.15791 & 0.93197 & -0.75755 \\
% 3 &  4\ 074 & 3.22361 & 3.16039 & 0.95586 & -0.97169 \\
% 4 &  5\ 856 & 3.22770 & 3.16240 & 0.98225 & -1.02418 \\
% 5 &  7\ 928 & 3.22724 & 3.16369 & 0.98459 & -1.02148 \\
% 6 & 10\ 290 & 3.22712 & 3.16401 & 0.98374 & -1.01846 \\
\midrule
$k = 1$ &  2\ 211 & 2.52594 & 2.47871 & 0.65728 & -0.81672 \\
$k = 2$ &  4\ 148 & 3.22841 & 3.16260 & 1.00571 & -1.03894 \\
$k = 3$ &  6\ 558 & 3.23184 & 3.16842 & 0.98822 & -1.02427 \\
$k = 4$ &  9\ 441 & 3.22714 & 3.16401 & 0.98431 & -1.01906 \\
$k = 5$ & 12\ 797 & 3.22759 & 3.16432 & 0.98578 & -1.02053 \\
$k = 6$ & 16\ 626 & 3.22757 & 3.16430 & 0.98580 & -1.02053 \\
\midrule
\multirow{2}{1.5cm}{\centering ref. \cite{featflow}}
& 167\ 232 & \markdigits{3.22}662 & \markdigits{3.16}351 & \markdigits{0.98}620 & \markdigits{-1.02}093 \\
& 667\ 264 & \markdigits{3.227}11 & \markdigits{3.164}26 & \markdigits{0.986}58 & \markdigits{-1.021}29 \\
% \midrule
% ref. \cite{schafer1996benchmark}&  & 3.22-3.24 &  & 0.99-1.01 & \\
\bottomrule
\end{tabular}
% END RECEIVE ORGTBL tab:2d
% \begin{comment}
% #+ORGTBL: SEND tab:2d orgtbl-to-latex :splice nil :skip 0
% |   k |  ndof |   max cD |   min cD |   max cL |    min cL |
% |-----+-------+----------+----------+----------+-----------|
% |   1 |  1380 | 2.736965 | 2.736965 | -0.15591 | -0.155909 |
% |   2 |  2582 | 3.207916 | 3.157910 |  0.93197 | -0.757550 |
% |   3 |  4074 | 3.223614 | 3.160391 |  0.95586 | -0.971686 |
% |   4 |  5856 | 3.227696 | 3.162395 |  0.98225 | -1.024181 |
% |   5 |  7928 | 3.227239 | 3.163690 |  0.98459 | -1.021481 |
% |   6 | 10290 | 3.227123 | 3.164014 |  0.98374 | -1.018457 |
% | ref |       | 3.227393 | 3.164263 |          |           |
% \end{comment}
  \end{center}
  \caption{Accuracy of the spatial discretization: results for different polynomial degrees.}
  \label{tab:2d}
\end{table}
We observe a rapid convergence for the $p$-refinement, i.e. the increase of the polynomial degree.
%, where the convergence in the drag values seems to be slightly better
Compared to the results from the literature \cite{featflow} the same order of accuracy is achieved with a lot less degrees of freedoms. 

\subsubsection{Numerical results: temporal discretization.} 
The results for the temporal discretization are shown in Table \ref{tab:time}. First of all, we observe that the second order IMEX scheme is already very accurate at its stability limit. But the method does not allow to choose larger time steps. This is in contrast to the alternative methods discussed here. For these methods we can consider much large time steps and observe a second order convergence. In this example the product decomposition method is more accurate than the modified fractional step method by one (time) level. 
We note that one major concern with this method is however, that splitting errors appear also if a stationary to the flow problem exists. This is not the case for additive decomposition methods or the proposed modified fractional step method.

\begin{table}[h!]
\small
  \begin{center}
% BEGIN RECEIVE ORGTBL tab:time
\begin{tabular}{lrrrlll}
\toprule
 & $1/\Delta t$ & $\max c_D$ & $\min c_D$ & $\max c_L$ & $\min c_L$ \\
\midrule
\multirow{3}{3cm}{\centering 2nd order IMEX} 
&$< $1 000    & \multicolumn{4}{c}{unstable}\\
&1 000    & 3.22754 & 3.16437 & 0.98486 & -1.01957 & \\
&32 000 & 3.22757 & 3.16431 & 0.98580 & -1.02053 \\
% \midrule
% 2nd order IMEX &1e-3    & 3.22754 & 3.16437 & 0.98486 & -1.01957 & \\
% 2nd order IMEX &5e-4    & 3.22754 & 3.16434 & 0.98522 & -1.01996 & \\
% 2nd order IMEX &2.5e-4  & 3.22755 & 3.16432 & 0.98551 & -1.02024 & \\
% 2nd order IMEX &1.25e-4 & 3.22756 & 3.16431 & 0.98567 & -1.02040 & \\
% \midrule
% 3rd order IMEX &2e-3 & 3.22747 & 3.16426 & 0.98365 & -1.01837\\
% 3rd order IMEX &1e-3 & 3.22719 & 3.16409 & 0.98346 & -1.01817\\
% 3rd order IMEX &5e-4 & 3.22714 & 3.16404 & 0.98360 & -1.01832\\
% % 3rd order IMEX &2.5e-4 & 3.22713 & 3.16402 &  & \\
\midrule
 \multirow{4}{3cm}{\centering operator-integration-\\factor splitting (BDF2)}
  &   125 & 3.38456 & 3.28279 & 1.16052 & -1.19939 \\
  &   250 & 3.25564 & 3.18643 & 1.01725 & -1.05276 \\
  &   500 & 3.23239 & 3.16836 & 0.98864 & -1.02378 \\
  & 1 000  & 3.22819 & 3.16394 & 0.98320 & -1.02013 \\
\midrule
 \multirow{3}{3cm}{\centering  modified fractional-\\ step-$\theta$-scheme}
 &   125 & 3.38272 & 3.29673 & 1.07667 & -1.12227 \\
 &   250 & 3.28713 & 3.21483 & 1.00937 & -1.04840 \\
 &   500 & 3.25036 & 3.18127 & 0.99172 & -1.02807 \\
 & 1 000 & 3.23656 & 3.17094 & 0.98721 & -1.02227 \\
% % \midrule
% % mod. Glowinksi &8e-3 & 3.47079 & 3.39370 & 1.07681 & -1.12582 \\
% % mod. Glowinksi &4e-3 & 3.34141 & 3.27385 & 1.01992 & -1.06153 \\
% % mod. Glowinksi &2e-3 & 3.26569 & 3.20140 & 0.99367 & -1.03122 \\
% % %mod. Glowinksi &1e-3 & 3.23652 & 3.17323 & 0.98490 & -1.02070 \\
% \midrule
% \texttt{FEATFLOW} \cite{featflow} & 2.5e-3   & \markdigits{3.227}39 & \markdigits{3.164}30 & \markdigits{0.986}37 & \markdigits{-1.021}08 \\
% \texttt{FEATFLOW} \cite{featflow} & 6.25e-4  & \markdigits{3.227}11 & \markdigits{3.164}26 & \markdigits{0.986}58 & \markdigits{-1.021}29 \\
% %(new ref glow:) & 1.25e-4 & 3.22825 & 3.16497 & 0.98585 & -1.02045 \\
% %ref. Table \ref{tab:2d}&  & 3.22712 & 3.16402 & 0.98374 & -1.01846 \\
% % \midrule
% % ref. \cite{schafer1996benchmark}& &  [3.22,3.24] &  & [0.99,1.01] & \\
\bottomrule
\end{tabular}
% END RECEIVE ORGTBL tab:time
  \end{center}
  \caption{Accuracy of operator-splitting time integration methods}
  \label{tab:time}
\end{table}

\subsection{A three-dimensional benchmark problem} \label{sec:numexbenchm}
Finally, we consider a three dimensional benchmark problem, the problem ``3D-3Z'' in \cite{schafer1996benchmark}. In contrast to the problem in section \ref{sec:numex2d} the inflow velocity is varied over time and the observed time interval is fixed. The maximal Reynolds number in this configuration is also $Re=100$ as in the previous section. The focus in this section is on the study of performance in the sense of computational effort over accuracy.

\subsubsection{Geometrical setup and boundary conditions.}
The geometrical setup is a generalization of the problem in the previous section. A cylindrical-shaped obstacle is places in a cuboid-shaped channel slightly above the vertical center:
\begin{equation}
  \Omega := [0,2.5] \! \times \! [0,0.41] \times \! [0,0.41] \setminus \{ \Vert (x_1,x_2) - (0.5,0.2) \Vert_2 \leq 0.05\}.
\end{equation}
 Boundary conditions are chosen as in the previous example except for a change to unsteady inflow boundary conditions:
\begin{equation*}
 u(0,y,z,t) = u_D(t) =  (9/4 \cdot \bar{u}(t)) \cdot 16 \cdot {y (d_y\!-\!y)}/{d_y^2} \cdot {z (d_z\!-\!z)}/{d_z^2} \cdot (1,0,0),
\end{equation*}
with $\bar{u}(t)$ the average inflow velocity which is time-dependent, $\bar{u}(t) = \sin(\pi t/8)$. The considered time interval is $[0,8s]$, the viscosity is set to $\nu = 10^{-3}$ s.t. the Reynolds number varies between $Re=0$ and $Re=100$. 

\subsubsection{Drag and Lift.} 
Again, the quantities of interest in this example are the (maximal and minimal) drag and lift forces $c_D$, $c_L$ that act on the disc. These are defined as
\begin{equation*}
 c_D := \frac{1}{\bar{u}_{\max}^2 r h } \int_{\Gamma_\circ} \left( \nu \frac{\partial u}{\partial \normal} - p \normal \right) \cdot e_x \ ds 
 \quad \quad  
 c_L := \frac{1}{\bar{u}_{\max}^2 r h} \int_{\Gamma_\circ} \left( \nu \frac{\partial u}{\partial \normal} - p \normal \right) \cdot e_y \ ds
\end{equation*}
with $\bar{u}_{\max} = \max_{t\in[0,8]} \bar{u}(t) = 1$, $r = 0.05$, $h = 0.41$ and $\Gamma_\circ$ the surface of the obstacle. 

\begin{figure}[h!]
  \begin{center}
    \includegraphics[width=0.55\textwidth]{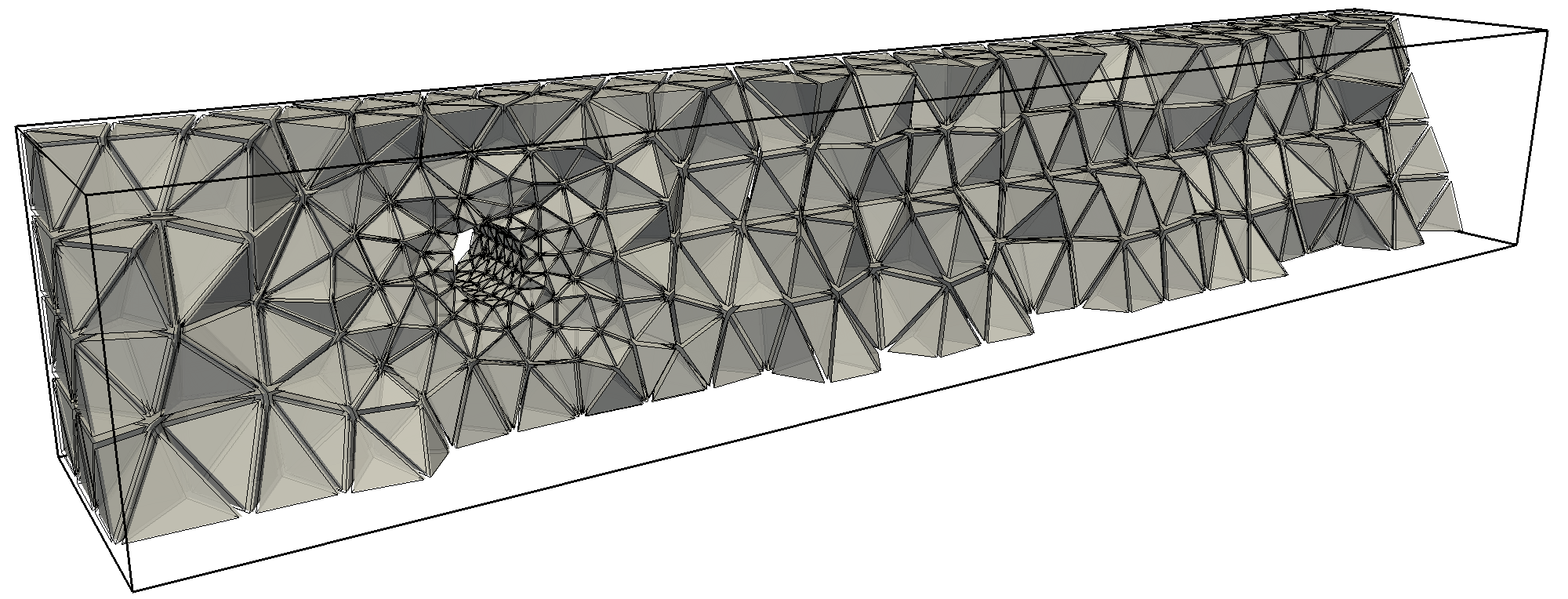}
    \includegraphics[trim=5cm 5cm 4cm 3cm, clip=true, width=0.44\textwidth]{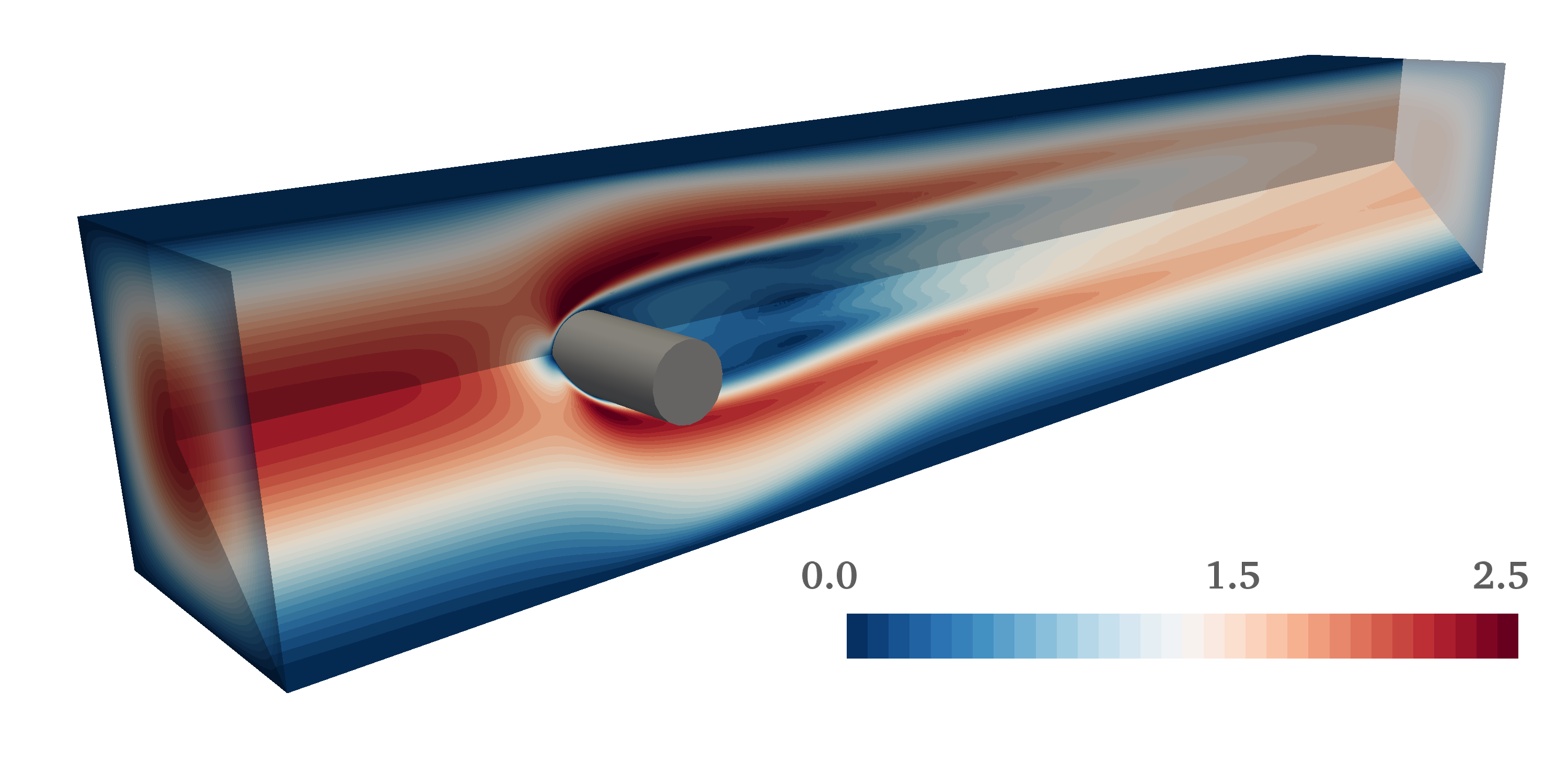}
  \end{center}
  \vspace*{-0.5cm}
  \caption{Used mesh (left) and solution at $t=0.4$ (right) for the benchmark problem in section \ref{sec:numexbenchm}.}
  \label{fig:sketch3d}
\end{figure}

% \begin{figure}[h!]
%   \begin{center}
%     \includegraphics[height=0.145\textwidth]{graphics/sol2d.png}
%     \hfill
%     \begin{overpic}[width=0.145\textwidth,angle=90]{graphics/scale.png}
%       \put(7,0){\small $0$}
%       \put(7,53){\small $2$}
%     \end{overpic}
%     \hfill
%     \hfill
%     \hfill
%     \includegraphics[height=0.145\textwidth]{graphics/blmesh_zoom.png}
%   \end{center}
%   \vspace*{-0.5cm}
%   \caption{Sketch of the mesh and the solution (color coding corresponding to velocity magnitude $\Vert \VelVar \Vert_2$) to the problem considered in section \ref{sec:numex2d} at a fixed time $t$ (left) and zoom-in on the boundary layer mesh (right).}
%   \label{fig:sketch2d}
% \end{figure}

\subsubsection{Numerical setup.}
We use an unstructured tetrahedral mesh consisting of 5922 elements, cf. Figure \ref{fig:sketch3d}. For the time discretization we use the modified fractional-step-$\theta$-scheme discussed in section \ref{opsplitting:glowinski}. To compensate for the increase in the spatial accuracy for increasing polynomial degree $k$, we adapt the number of time steps accordingly. The computations were carried out on a shared-memory computer with $24$ cores. We comment on details of the computing times in remark \ref{rem:timings}.

\begin{table}[h!]
\small
  \begin{center}
% BEGIN RECEIVE ORGTBL tab:time
\begin{tabular}{crrrrlr}
\toprule
 & \#ndof [K] & $\max c_D$ & $\max c_L$ & $\min c_L$ &$\Delta t$ [s] & comp. time [s]\\
% \midrule
% $k=3$ & 3.33390 & 0.00282 & -0.013165 & 194 & 0.0025 & 2 008 $\times$ 24\\
% $k=4$ & 3.29844 & 0.00278 & -0.010902 & 290 & 0.001  & 7 846 $\times$ 24\\
% $k=5$ & 3.29921 & 0.00278 & -0.010889 & 407 & 0.0025 & 7 577 $\times$ 24\\
% $k=6$ & 3.29863 & 0.00278 & -0.011033 & 542 & 0.0025 & 14 349 $\times$ 24\\
% \midrule
% older:
% $k=2$ &  169 & 3.43053 & 0.00262 & -0.016289 & 0.01   &   418 $\times$ 24\\
% $k=3$ &  343 & 3.29330 & 0.00277 & -0.011099 & 0.0080 & 1 770 $\times$ 24\\
% $k=4$ &  595 & 3.29856 & 0.00278 & -0.010762 & 0.0050 &${}^\ast$ 2 604 $\times$ 24\\
% $k=5$ &  939 & 3.29795 & 0.00278 & -0.011053 & 0.0025 & 8 745 $\times$ 24\\
\midrule
% recent:
$k=2$ &  169 & 3.43046 & 0.00262 & -0.016289 & 0.0080 &    492 $\times$ 24\\
$k=3$ &  343 & 3.29331 & 0.00277 & -0.011099 & 0.0080 &    964 $\times$ 24\\
$k=4$ &  595 & 3.29853 & 0.00278 & -0.010762 & 0.0040 &  3 087 $\times$ 24\\
$k=5$ &  939 & 3.29798 & 0.00278 & -0.011054 & 0.0040 &  6 670 $\times$ 24\\
\midrule
% ref. \cite{bayraktar2012benchmark} (Lv1) & 3.2207\hphantom{0}  & 0.0027\hphantom{0}  & -0.009500 &    199 & 0.01    &   3 220 $\times$ \hphantom{0}2\\
% ref. \cite{bayraktar2012benchmark} (Lv2) & 3.2877\hphantom{0}   & 0.0028\hphantom{0} & -0.010892 &  1 482 & 0.01    &  17 300 $\times$ \hphantom{0}4\\
\multirow{2}{2cm}{\centering ref. \cite{bayraktar2012benchmark}}  
& 11 432 & 3.2963\hphantom{0}   & 0.0028\hphantom{0} & -0.010992  & 0.01    &  35 550 $\times$ 24\\
& 89 760 & 3.2978\hphantom{0}   & 0.0028\hphantom{0} & -0.010999 & 0.005   & 214 473 $\times$ 48\\
\midrule
ref. \cite{john2006efficiency}& 7 036 & 3.2968\hphantom{0} & & -0.011\hphantom{000} &  & \\
\midrule
ref. \cite{schafer1996benchmark}& & [3.2,3.3] & [0.002,0.004] &  &  \\
\bottomrule
\end{tabular}
% END RECEIVE ORGTBL tab:time
  \end{center}
  \caption{Numerical results for the benchmark problem ``3D3Z'' in \cite{schafer1996benchmark}. Results obtained with different polynomial degrees and reference values. }
  \label{tab:bench3d}
\end{table}

\subsubsection{Numerical results.}
In Table \ref{tab:bench3d} the results obtained are compared with the literature in terms of accuracy \emph{ and computing time}. We observe that we can achieve the same level of accuracy as the results in \cite{bayraktar2012benchmark} (and \cite{john2006efficiency}) with a computing time which is dramatically smaller. We note that the computing time per degree of freedom is actually worse than in \cite{bayraktar2012benchmark}. Nevertheless, the same accuracy is achieved with much less degrees of freedoms using the high order method ($k>2$), s.t. our computations exceed the performance results in the literature.
In the study \cite{bayraktar2012benchmark} one of the conclusions is that their third order method ($Q_2/P_1^{\text{disc}}$) is much more efficient compared to lower order methods. We extend this conclusion in the sense that the use of even higher order discretizations, i.e. $k>2$, increases efficiency even further. Moreover, high order discretizations can be implemented efficiently.
One important component for the efficient handling of the Navier-Stokes equations with our high order discretization is the time integration using operator-splitting.

\begin{remark}[Computation times]\label{rem:timings}
We remark on the computation for the case $k=5$. In that computation approximately $65\%$ of the computing time has been spend on the solution of linear systems for Stokes-type problems, $30\%$ on convection operator evaluations and $5\%$ on the setup and remaining operations. To solve the implicit convection problem (Step 2 in \eqref{eq:modglow}) an average of $20$ iterations has been applied.
\end{remark}

\begin{remark}[Linear systems]\label{rem:precond}
In the test cases in this section we only applied direct solvers which is possible due to the (comparably) small size of the arising linear systems. For problems with increasing complexity efficient linear solvers are mandatory. The development of suitable preconditioners of the Stokes problem is based on efficient preconditioning of the bilinear form $\bilinearform{A}$. For the scalar problem we could show poly-logarithmic bounds in $k$ for the condition number of standard p-version domain decomposition preconditioners, cf. \cite{schoberl2013domain}. We plan to investigate suitable generalizations of this preconditioner for the Stokes problem in the future.
\end{remark}

\section{Conclusion} \label{sec:concl}
We presented and discussed a combined DG/HDG discretization tailored for efficiency. We summarize the core components. We split the Navier-Stokes problem into linear Stokes-type problems and hyperbolic transport problems by means of operator-splitting time integration. For the Stokes-type problem we use an $H(\divergence)$-conforming Hybrid DG formulation with a new modification: the \emph{projected jumps} formulation. 
The Hybrid DG formulation facilitates the efficient solution of linear systems compared to other DG methods. The \emph{projected jumps} formulation improves its efficiency even more. 
For the hyperbolic transport problem we apply a standard DG formulation. In numerical test cases we demonstrated the performance of the method.
\section*{Acknowledgements}
The authors greatly appreciate the contribution of Philip Lederer related to the calculation of LBB-constants, cf. remark \ref{rem:probust}.

\bibliographystyle{wileyj}
\bibliography{hdiv_hdg_nse}

\end{document}